\input amssym.def
\input amssym.tex


\def\item#1{\vskip1.3pt\hang\textindent {\rm #1}}


\tolerance=300
\pretolerance=200
\hfuzz=1pt
\vfuzz=1pt


\hoffset=0.6in
\voffset=0.8in

\hsize=5.8 true in 
\vsize=8.5 true in
\parindent=25pt
\mathsurround=1pt
\parskip=1pt plus .25pt minus .25pt
\normallineskiplimit=.99pt

\countdef\revised=100
\mathchardef\emptyset="001F 
\chardef\ss="19
\def\3{\ss}
\def\anf{$\lower1.2ex\hbox{"}$}
\def\frac#1#2{{#1 \over #2}}
\def\>{>\!\!>}
\def\<{<\!\!<}

\def\into{\hookrightarrow}
\def\ssarr{\hbox to 30pt{\rightarrowfill}}
\def\sarr{\hbox to 40pt{\rightarrowfill}}
\def\arr{\hbox to 60pt{\rightarrowfill}}
\def\larr{\hbox to 60pt{\leftarrowfill}}
\def\Arr{\hbox to 80pt{\rightarrowfill}}

\def\ssmapright#1{\smash{\mathop{\ssarr}\limits_{#1}}}

{}

\def\Ad{\mathop{\rm Ad}\nolimits}

\def\conv{\mathop{\rm conv}\nolimits}

\def\End{\mathop{\rm End}\nolimits}

\def\Gl{\mathop{\rm Gl}\nolimits}

\def\Hom{\mathop{\rm Hom}\nolimits}%
\def\im{\mathop{\rm im}\nolimits}
\def\Im{\mathop{\rm Im}\nolimits}

\def\Int{\mathop{\rm int}\nolimits}

\def\Re{\mathop{\rm Re}\nolimits}

\def\Sl{\mathop{\rm Sl}\nolimits}
\def\SO{\mathop{\rm SO}\nolimits}

\def\sup{\mathop{\rm sup}\nolimits}
\def\supp{\mathop{\rm supp}\nolimits}


\def\0{{\bf 0}}
\def\1{{\bf 1}}

\def\a{{\frak a}}

\def\b{{\frak b}}

\def\e{{\frak e}}

\def\g{{\frak g}}

\def\h{{\frak h}}

\def\k{{\frak k}}

\def\n{{\frak n}}

\def\p{{\frak p}}
\def\q{{\frak q}}

\def\C{{\Bbb C}} 
\def\D{{\Bbb D}}

\def\N{{\Bbb N}} 
 
\def\P{{\Bbb P}} 
 
\def\R{{\Bbb R}}

\def\:{\colon}  
\def\.{{\cdot}}
\def\|{\Vert}
\def\bsk{\bigskip}

\def\giantskip{\vskip2\bigskipamount}
\def\gsk{\giantskip}
\def \la {\langle}
\def\msk{\medskip}
\def \ra {\rangle}
\def \res {\!\mid\!\!}

\def\ssk{\smallskip}

\def\bbr{\bigbreak}
\def\giantbreak{\par \ifdim\lastskip<2\bigskipamount \removelastskip
         \penalty-400 \giantskip\fi}

\def\nin{\noindent}
\def\cen{\centerline}
\def\pagebreak{\vskip 0pt plus 0.0001fil\break}
\def\linebreak{\break}

\def\hat{\widehat}

\def\eps{\varepsilon}
\def\epsilon{\varepsilon}

\def\nin{\noindent}
\def\oline{\overline}

\def\pder#1,#2,#3 { {\partial #1 \over \partial #2}(#3)}
\def\pde#1,#2 { {\partial #1 \over \partial #2}}
\def\phi{\varphi}


\def\subeq{\subseteq}
\def\supeq{\supseteq}

\def\tilde{\widetilde}

\def\up{{\uparrow}}

\font\eightrm=cmr8


\font\bfone=cmbx10 scaled\magstep1 
\font\bftwo=cmbx10 scaled\magstep2 

\def\qed{{\unskip\nobreak\hfil\penalty50\hskip .001pt \hbox{}\nobreak\hfil
          \vrule height 1.2ex width 1.1ex depth -.1ex
           \parfillskip=0pt\finalhyphendemerits=0\medbreak}\rm}


\def\Lemma #1. {\bigbreak\vskip-\parskip\noindent{\bf Lemma #1.}\quad\it}

\def\Sublemma #1. {\bigbreak\vskip-\parskip\noindent{\bf Sublemma #1.}\quad\it}

\def\Proposition #1. {\bigbreak\vskip-\parskip\noindent{\bf Proposition #1.}
\quad\it}

\def\Corollary #1. {\bigbreak\vskip-\parskip\nin{\bf Corollary #1.}
\quad\it}

\def\Theorem #1. {\bigbreak\vskip-\parskip\noindent{\bf Theorem #1.}
\quad\it}

\def\Definition #1. {\rm\bigbreak\vskip-\parskip\noindent{\bf Definition #1.}
\quad}

\def\Remark #1. {\rm\bigbreak\vskip-\parskip\noindent{\bf Remark #1.}\quad}

\def\Example #1. {\rm\bigbreak\vskip-\parskip\noindent{\bf Example #1.}\quad}

\def\Problems #1. {\bigbreak\vskip-\parskip\noindent{\bf Problems #1.}\quad}
\def\Problem #1. {\bigbreak\vskip-\parskip\noindent{\bf Problems #1.}\quad}

\def\Conjecture #1. {\bigbreak\vskip-\parskip\noindent{\bf Conjecture #1.}\quad}

\def\Proof#1.{\rm\par\ifdim\lastskip<\bigskipamount\removelastskip\fi\smallskip
            \noindent {\bf Proof.}\quad}

\def\Axiom #1. {\bigbreak\vskip-\parskip\noindent{\bf Axiom #1.}\quad\it}

\def\Satz #1. {\bigbreak\vskip-\parskip\noindent{\bf Satz #1.}\quad\it}

\def\Korollar #1. {\bbr\vskip-\parskip\nin{\bf Korollar #1.} \quad\it}

\def\Bemerkung #1. {\rm\bigbreak\vskip-\parskip\noindent{\bf Bemerkung #1.}
\quad}

\def\Beispiel #1. {\rm\bigbreak\vskip-\parskip\noindent{\bf Beispiel #1.}\quad}
\def\Aufgabe #1. {\rm\bigbreak\vskip-\parskip\noindent{\bf Aufgabe #1.}\quad}

\def\Beweis#1. {\rm\par\ifdim\lastskip<\bigskipamount\removelastskip\fi
           \smallskip\noindent {\bf Beweis.}\quad}

\nopagenumbers

\def\date{\ifcase\month\or January\or February \or March\or April\or May
\or June\or July\or August\or September\or October\or November
\or December\fi\space\number\day, \number\year}

\def\title{Title ??}
\def\author{Author ??}

\def\thanks#1{\footnote*{\eightrm#1}}

\def\rightheadline{\hfil{\eightrm\title}\hfil\tenbf\folio}
\def\leftheadline{\tenbf\folio\hfil{\eightrm\author}\hfil}
\headline={\vbox{\line{\ifodd\pageno\rightheadline\else\leftheadline\fi}}}

\def\firstheadline{}
\def\firstfootline{\cen{\rm\folio}}

\def\seite #1 {\pageno #1
               \headline={\ifnum\pageno=#1 \firstheadline
               \else\ifodd\pageno\rightheadline\else\leftheadline\fi\fi}
               \footline={\ifnum\pageno=#1 \firstfootline\else{}\fi}}

\newdimen\dimenone
 \def\checkleftspace#1#2#3#4{
 \dimenone=\pagetotal
 \advance\dimenone by -\pageshrink   
 \ifdim\dimenone>\pagegoal          
   \else\dimenone=\pagetotal
        \advance\dimenone by \pagestretch
        \ifdim\dimenone<\pagegoal
          \dimenone=\pagetotal
          \advance\dimenone by#1         
          \setbox0=\vbox{#2\parskip=0pt                
                     \hyphenpenalty=10000
                     \rightskip=0pt plus 5em
                     \noindent#3 \vskip#4}    
        \advance\dimenone by\ht0
        \advance\dimenone by 3\baselineskip   
        \ifdim\dimenone>\pagegoal\vfill\eject\fi
          \else\eject\fi\fi}


\def\subheadline #1{\nin\bigbreak\vskip-\lastskip
      \checkleftspace{0.7cm}{\bf}{#1}{\medskipamount}
          \indent\vskip0.7cm\centerline{\bf #1}\medskip}

\def\sectionheadline #1{\bigbreak\vskip-\lastskip
      \checkleftspace{1.1cm}{\bf}{#1}{\bigskipamount}
         \vbox{\vskip1.1cm}\cen{\bfone #1}\bsk}

\def\lsectionheadline #1 #2{\bigbreak\vskip-\lastskip
      \checkleftspace{1.1cm}{\bf}{#1}{\bigskipamount}
         \vbox{\vskip1.1cm}\cen{\bfone #1}\msk \cen{\bfone #2}\bsk}

\def\lchapterheadline #1 #2{\bigbreak\vskip-\lastskip\indent\vskip3cm
                       \cen{\bftwo #1} \msk \cen{\bftwo #2} \gsk}
\def\llsectionheadline #1 #2 #3{\bigbreak\vskip-\lastskip\indent\vskip1.8cm
\cen{\bfone #1} \msk \cen{\bfone #2} \msk \cen{\bfone #3} \nobreak\bsk\nobreak}


\newtoks\literat
\def\[#1 #2\par{\literat={#2\unskip.}%
\hbox{\vtop{\hsize=.15\hsize\nin [#1]\hfill}
\vtop{\hsize=.82\hsize\nin\the\literat}}\par
\vskip.3\baselineskip}

\mathchardef\emptyset="001F 
\def\address{Author: \tt$\backslash$def$\backslash$address$\{$??$\}$}

\def\firstpage{\nin
{\obeylines \parindent 0pt }
\vskip2cm
\centerline {\bfone \title}
\gsk
\centerline{\bf\author}

\vskip1.5cm \rm}

\def\addresstwo{}

\def\dlastpage{\par\vbox{\vskip1cm\nin
\line{
\vtop{\hsize=.5\hsize{\parindent=0pt\baselineskip=10pt\nin\address}}
\quad 
\vtop{\hsize=.42\hsize\nin{\parindent=0pt
\baselineskip=10pt\addresstwo}}
\hfill} }}


\def\b{\mathop{\bf b}\nolimits}

\def\c{\mathop{\bf c}\nolimits}

\def\cH{\mathop{{\bf c}_{\hbox{\fiverm G/H}}}\nolimits}

\def\e{\mathop{\bf e}\nolimits}

\def\hfH{\hbox{\fiverm H}}

\def\hfK{\hbox{\fiverm K}}

\def\ssmapright#1{\smash{\mathop{\ssarr}\limits^{#1}}}

\input epsf.tex

\def\bs{\backslash} 
 
\def\addots{\mathinner{\mkern1mu\raise1pt\vbox{\kern7pt\hbox{.}}\mkern2mu 
\raise4pt\hbox{.}\mkern2mu\raise7pt\hbox{.}\mkern1mu}} 
 
\pageno=1 
\def\up#1{\leavevmode \raise.16ex\hbox{#1}} 
 at 8truept 
 at 8truept 
 at 12truept 
\chardef\ss="19 
\def\3{\ss} 
 
\def\firstpage{\nin 
{\obeylines \parindent 0pt } 
\vskip2cm 
\centerline {\bfone \title} 
\ssk 
\centerline {\bfone \titletwo} 
\gsk 
\centerline{\bf\author} 
\vskip1.5cm \rm}

\def\title{Holomorphic $H$-spherical distribution vectors} 
\def\titletwo{in principial series representations} 
 
\def\author{Simon Gindikin, Bernhard Kr\"otz and 
Gestur \'Olafsson} 
 
\footnote{}{SG was supported in part  by  NSF-grant DMS-0070816} 
\footnote{}{BK was supported in part  by NSF-grant DMS-0097314} 
\footnote{}{G\'O was supported in part by NSF-grant  DMS-0070607 and DMS-0139783}

\def\address 
{Simon Gindikin 
 
Department of Mathematics 
 
Rutgers University 
 
New Brunswick, NJ 08903 
 
USA 
 
{\tt gindikin@math.rutgers.edu} 
 
\bsk 
\bsk 
 
Gestur \'Olafsson 
 
Louisiana State University 
 
Department of Mathematics 
 
Baton Rouge, LA 70803 
 
USA 
 
{\tt olafsson@math.lsu.edu} 
 
} 
 
\def\addresstwo 
{Bernhard Kr\"otz 
 
Department of Mathematics 
University of Oregon  
Eugene Or 97403-1221  
USA  
 
{\tt kroetz@math.uoregon.edu} 
}

\firstpage

\subheadline{Abstract} 
\noindent 
Let $G/H$ be a semisimple symmetric space. 
The main tool to embed a principal series 
representation of  $G$ into $L^2(G/H)$ are 
the $H$-invariant distribution vectors. 
If $G/H$ is a non-compactly 
causal symmetric space, then $G/H$ can be realized as a 
boundary component of the complex crown $\Xi$. 
In this article we construct a minimal $G$-invariant 
subdomain $\Xi_H$ of $\Xi$ with 
$G/H$ as Shilov boundary. Let $\pi$ be 
a spherical principal series representation 
of $G$. We show 
that the space of $H$-invariant 
distribution vectors of $\pi$, which admit a 
holomorphic extension to $\Xi_H$,  
is one dimensional. Furthermore we give a spectral definition of 
a Hardy space corresponding to those distribution 
vectors. In particular we achieve a geometric 
realization of a multiplicity free subspace of $L^2(G/H)_{\rm mc}$ 
in a space of holomorphic functions. 
 
\sectionheadline{Introduction} 
 
\noindent 
Holomorphic extensions and boundary value maps have been 
valuable tools to solve problems in representation theory 
and harmonic analysis on {\it real} symmetric 
spaces. Two of the best known constructions  
are Hardy spaces with their boundary value maps and  
Cauchy-Szeg\"o-kernels, 
and Fock space constructions with their corresponding 
Segal-Barmann transform. It is in this flavour that we establish a 
correspondence between eigenfunctions on a Riemannian symmetric spaces $X=G/K$ 
and  a  non-compactly causal (NCC) symmetric spaces $Y=G/H$ in  this paper. 
In particular 
we, via analytic continuation, relate a {\it spherical function} $\phi_\lambda$ on 
$G/K$ to a {\it holomorphic} $H$-invariant distribution on $G/H$.

\par Let us explain our results in more detail. On the geometric level we construct a certain minimal 
$G$-invariant Stein domain $\Xi_H\subeq X_\C=G_\C/K_\C$ with 
the following properties: The Riemannian symmetric space $X$ is embedded into $\Xi_H$ as 
a totally real submanifold and the affine non-compactly causal space $Y$ is isomorphic to the 
distinguished (Shilov) boundary of $\Xi_H$. The details 
of this construction are carried out in Section 1. 
 
\par The minimal tube $\Xi_H$ is a subdomain of the complex 
crown $\Xi\subeq X_\C$ of $X$ -- an object first introduced 
in [AG90] which became subject of intense study over the last 
few years. A consequence is that all $\D(X)$-eigenfunctions 
on $X$ extend holomorphically to $\Xi_H$ [KS01b]. Another key fact 
is that $\D(X)\simeq \D (Y)$. Thus by taking limits on the boundary $Y$ we  obtain a realization 
of the $\D(X)$-eigenfunctions on $X$ as $\D(Y)$-eigenfunctions on $Y$. 
Conversely, eigenfunctions on $Y$ which holomorphically extend to $\Xi_H$ 
yield by restriction eigenfunctions on $X$.

\par It seems to us that the above mentioned transition between eigenfunctions 
on $X$ and $Y$ is most efficiently described using the 
techniques from representation theory.  To fix the notation let $(\pi,{\cal H})$ 
denote an admissible Hilbert representation of $G$ with 
finite length. We write ${\cal H}^K$ for the space of $K$-fixed 
vectors and $({\cal H}^{-\infty})^H$ for the space of $H$-fixed 
distribution vectors of $\pi$. Using the method of analytic 
continuation of representations as developed in [KS01a] we 
establish a bijection 
$${\cal H}^K\ssmapright{\simeq} ({\cal H}^{-\infty})_{\rm hol}^H , \ \ 
v_{\hbox{\fiverm K}}\mapsto v_{\hbox{\fiverm H}}$$ 
where $({\cal H}^{-\infty})_{\rm hol}^H\subeq ({\cal H}^{-\infty})^H$ 
denotes the subspace characterized through the property that 
associated matrix coefficients on $Y$ extend holomorphically 
to $\Xi_H$ (cf.\ Theorem 2.1.3, Theorem 2.2.4). This bijection and various 
ramifications are the subject proper of Section 2. 
 
\par In Section 3 we give an application 
of our theory towards the geometric realization of the most-continuous spectrum 
$L^2(Y)_{\rm mc}$ of $L^2(Y)$. First progress in this direction was achieved in [GK\'O01]. 
There, for the cases where $\Xi=\Xi_H$,  we defined a Hardy space 
${\cal H}^2(\Xi)$ on $\Xi$ and showed that there is an isometric boundary value mapping 
realizing ${\cal H}^2(\Xi)$ as a 
multiplicity one subspace of $L^2(Y)_{\rm mc}$ of full spectrum. 
It was an open problem how to define 
Hardy spaces for general NCC symmetric spaces $Y$ and to determine  the Plancherel 
measure explicitely. We solve this problem by giving a spectral definition of the 
Hardy space, i.e., we take the conjectured Plancherel measure and define a Hilbert 
space of holomorphic functions  ${\cal H}^2(\Xi_H)$ on $\Xi_H$. The identification 
of ${\cal H}^2(\Xi_H)$ as  a Hardy space then follows by establishing 
an isometric boundary value mapping $b\: {\cal H}^2(\Xi_H)\into L^2(G/H)_{\rm mc}$. 
In particular we achieve a geometric realization of a multiplicity free subspace of $L^2(Y)_{\rm mc}$ 
in holomorphic functions. 
 
\msk It is our pleasure to thank the referee for his very careful  work. He pointed out many  
inaccuracies and made useful remarks on the presentation of the paper.

\sectionheadline{1. Complex crowns and the domains $\Xi_H$} 
 
\noindent 
The purpose of this section is to give the geometric preliminaries  
of the analytical constructions to come. Our two main  
players are a Riemannian symmetric space $G/K$ on the one hand side  
and on the other hand a non-compactly causal symmetric space $G/H$.  
The two symmetric spaces $G/K$ and $G/H$ are ``connected'' through  
a complex $G$-invariant domain $\Xi_H\subeq G_\C/K_\C$ in the  
following way: $G/K$ is a totally real submanifold and  
$G/H$ constitutes the distinguished (Shilov) boundary of $\Xi_H$.  
The domain $\Xi_H$ constructed in this section is an appropriate  
subdomain of the complex crown $\Xi$ of the Riemannian symmetric space 
$G/K$. 
 
\par This section is organized as follows. We start by briefly recalling  
the defintion and some standard features of non-compactly causal  
symmetric spaces . Then we switch to complex crowns $\Xi$  
and summarize the main results of [GK02a] on how 
to realize $G/H$ in the distinguished boundary of $\Xi$. Finally  
we give the construction of the domain $\Xi_H$.

\subheadline{1.1. Non-compactly causal symmetric spaces (NCC)} 
 
\noindent 
In this subsection we recall some facts on non-compactly causal symmetric spaces. 
The material is standard and can be found in the monograph 
[H\'O96]. 
 
\ssk  
Let $G$ be a connected semisimple  Lie group and $\g$ be its Lie algebra. 
Denote by $\g_\C=\g\otimes_\R\C$ the complexification of $\g$. 
If $\h$ is a subalgebra of $\g$,  then we denote by $\h_\C$ the  complex subalgebra 
of $\g_\C$ generated by $\h$. 
We assume that  $G$ is contained 
in a complex group $G_\C$ with Lie algebra $\g_\C$. 
 
\par If $\sigma\:G \to G$ is 
an involution, then, by abuse of notation, we use the same letter 
for the derived involution on the Lie algebra $\g$ and its complex linear  
extension to  $\g_\C$. 
\par Let $\theta\:G\to G$ be a Cartan involution and denote by $K<G$ the corresponding 
maximal compact subgroup. Let $\k=\{X\in\g\: \theta(X)=X\}$ 
and $\p =\{X\in \g\:  \theta (X)=-X\}$. Then 
$\k$ is the Lie algebra of $K$. 
 
\par In the sequel we let $\tau$ denote an involution on $G$ which we may assume to commute 
with $\theta$. Let $G^\tau :=\{g\in G\: \tau (g)=g\}$ and 
let $H$ be an open subgroup of $G^\tau$. Then $G/H$ is called  a {\it symmetric 
spaces}.  On the Lie algebra level 
$\tau$ induces a splitting $\g=\h+\q$ with $\h$ the $+1$ and $\q$ the $-1$-eigenspace of 
$\tau$. Notice that $\h$ is the Lie 
algebra of $H$. The pair $(\g,\h)$ is called a {\it symmetric pair}. We have, as $\theta$ and $\tau$ commute: 
$$\eqalign{ \g &=\k+ \p\cr 
&=\h+ \q\cr 
&=\k\cap \h +\k\cap \q + \p\cap \h + \p\cap \q 
} 
$$ 
\par 
The symmetric pair 
$(\g ,\h)$ is called {\it irreducible} if the only 
$\tau$-invariant ideals in $\g$ are the trivial ones, 
$\{0\}$ and $\g$. In this case either $\g$ is 
simple or $\g \simeq \h\oplus \h$, with $\h$ simple, and $\tau (X,Y)= 
(Y,X)$ the flip. We say that the symmetric space $G/H$ is {\it irreducible} 
if the corresponding symmetric pair $(\g ,\h)$ is irreducible. 
 
\par Let $\emptyset \not= C\subseteq \g$ be an open subset of $\g$. Then 
$C$ is said to be {\it hyperbolic} if for all $X\in C$ the map ${\rm ad}(X):\g \to \g$ is 
semisimple with real eigenvalues. 
 
\Definition 1.1.1. {\bf (NCC)} Assume that $G/H$ is an irreducible symmetric space. Then 
the following two conditions are equivalent: 
 
\item{(a)} There exists a non-empty $H$-invariant open hyperbolic convex cone $C 
\subseteq \q$ which contains  no affine lines; 
 
\item{(b)} There exists an element $T_0\in \q\cap \p$, $T_0\not= 0$, 
which is fixed by $H\cap K$. 
 
If one of those equivalent conditions are satisfied, then $G/H$ is called 
{\it non-compactly causal}, or {\it NCC} for short. 
\qed 
 
\Remark 1.1.2. (a) The element $T_0$ in Definition 1.1.1 is 
unique up to multiplication by scalar. We can normalize 
$T_0$ such that ${\rm ad}(T_0)$ has spectrum $\{0,1,-1\}$. The eigenspace 
corresponding to $0$ is exactly $\g^{\theta\tau}=\k\cap\h +\p\cap \q$. 
\par\nin (b) If $G/H$ is NCC and $\a\subset \p\cap \q$ is maximal 
abelian, then $T_0\in \a$ by (a). Hence, again by (a), it follows 
that $\a$ is also maximal abelian in $\p$ and in $\q$. 
\par\nin (c) Let $T_0$ be as above. Then the interior of the convex hull of $\R^+\Ad (H)T_0$ 
is a minimal, $H$-invariant open hyperbolic convex cone in $\q$. 
\par\nin (d) All the NCC pairs $(\g, \h)$ are classified and we refer 
to [H\'O96, Th.\ 3.2.8] for the complete list.\qed

\subheadline{1.2. The complex crown of a Riemannian symmetric space} 
 
\noindent 
The NCC spaces are exactly the affine symmetric spaces that can be 
realized as a symmetric subspace in the distinguished boundary 
of the complex crown $\Xi$ of the Riemannian symmetric space 
$G/K$. We will therefore recall some basic facts about 
$\Xi$. We refer to [GK02a] and [GK02b] as a standard source.

Let the notation be as in Subsection 1.1. Let $\a\subeq \p$ be a maximal 
abelian subalgebra. For $\alpha \in \a^*$ let 
$\g^\alpha=\{ X\in \g\: (\forall H\in \a)\ [H,X]=\alpha(H)X\}$ and 
let $\Sigma:=\{\alpha\in \a^*\:\alpha\not=0,\g^\alpha\not=\{0\}\}$ be 
the corresponding set of restricted roots.

\par Following [AG90] we define a bounded convex subset of $\a$ by 
$$\Omega=\{ X\in \a\: (\forall \alpha\in \Sigma)\ |\alpha(X)|<{\pi\over 2}\}\ .$$ 
Denote by $K_\C$ the analytic subgroup of $G_\C$ with Lie algebra $\k_\C$. 
Then we define a $G-K_\C$ double coset domain in $G_\C$ by 
 
$$\tilde \Xi=G\exp(i\Omega)K_\C$$ 
and recall that $\tilde\Xi$ is open in $G_\C$ [KS01a]. In particular the domain 
$$\Xi=\tilde\Xi/K_\C$$ 
is an open $G$-invariant subset of  $G_\C/K_\C$ containing $G/K$ as a totally 
real submanifold. We refer to $\Xi$ as the {\it complex crown} of the 
Riemannian symmetric space $G/K$ (cf.\ [AG90]). Observe that the definition 
of $\Xi$ and $\tilde \Xi$ is independent of the choice of $\a\subeq \p$. 
For a subset $\omega\subeq \a$ we define a tube domain in $A_\C=\exp(\a_\C)$ 
by 
$$T( \omega) =A\exp(i\omega)$$ 
and notice that $T (2\Omega)$ is biholomorphic to $\a+i 2\Omega$ via the 
exponential map. 
 
\msk Fix a positive system $\Sigma^+$ of $\Sigma$ and define a subalgebra 
$\n$ of $\g$ by 
$$\n=\bigoplus_{\alpha\in \Sigma^+} \g^\alpha\ .$$ 
Write $N_\C$ for the analytic subgroup of $G_\C$ with Lie algebra $\n_\C$. 
Then it follows from [GK02b] that 
$$\tilde\Xi\subeq N_\C T(\Omega) K_\C\leqno(1.2.1)$$ 
or even more precisely 
$$\tilde\Xi=\left[\bigcap_{g\in G} g  N_\C T(\Omega) K_\C \right]_0\leqno (1.2.2)$$ 
where the subscript ${}_0$ denotes the connected component of $[\cdot]$ containing $G$. 
 
\subheadline{1.3. The distinguished boundary of $\Xi$} 
 
\noindent 
Write $\oline \Xi$ and $\partial \Xi$ for the closure respectively the boundary  of $\Xi$ in 
$G_\C/K_\C$. The {\it distinguished boundary} of $\Xi$ is a certain finite 
union of $G$-orbits in $\partial \Xi$ which features many properties of 
a Shilov boundary. It was introduced and investigated in [GK02a] 
and the objective of this subsection is to recall its definition and basic 
properties.

\par  Let 
${\cal W}=N_K(\a)/ Z_K(\a)$ be the Weyl group of $\a$ in $G$. 
Write $\oline\Omega$ for the closure of $\Omega$ and notice that 
$\oline\Omega$ is a ${\cal W}$-invariant compact convex set. Denote by 
$\partial_e\Omega$ the set of extreme points of $\oline \Omega$. Then 
there exists $X_1,\ldots ,X_n\in\partial_e\Omega$ such 
that 
$$\partial_e\Omega={\cal W}(X_1)\amalg\ldots\amalg{\cal W}(X_n)\ .$$ 
Define the {\it distinguished boundary}  of $\Xi$ in $G_\C/K_\C$ 
by 
$$\partial_d\Xi\:=G\exp(i\partial_e\Omega)K_\C/K_\C\ .$$ 
We refer to [GK02a] for detailed information about $\partial_d\Xi$ 
and recall here only the facts that we need. 
 
\par For $1\leq j\leq n$  let 
$z_j\:=\exp(iX_j)K_\C \in \partial_d\Xi$. If $G_\C$ is not simply connected it 
can happen that $G(z_j)=G(z_k)$ for some $j\neq k$. But after relabelling 
the $z_j$ we can assume that there is an $m\leq n$ such that 
$G(z_j)\not=G(z_k)$ for $1\le j,k\le m$, $j\not= k$, and 
$$\partial_d\Xi= G(z_1)\amalg\ldots\amalg G(z_m)\ .$$ 
Denote by $H_j$  the isotropy subgroup of $G$ in $z_j$. Then as $G$-spaces we have 
 
$$\partial_d\Xi= G/H_1\amalg\ldots \amalg G/ H_m \ .$$ 
As a consequence of the complete classification of $\partial_d\Xi$ in 
[GK02a] we obtain the following fact.

\Proposition 1.3.1. For the distinguished boundary  
$\partial_d\Xi$ of $\Xi$  the following assertions hold: 
\item{(i)} If one of the boundary components 
$G/H_j$ of $\partial_d\Xi$ is a symmetric space, then it is a 
non-compactly causal symmetric space.  
\item{(ii)} Every 
non-compactly causal symmetric space of the form $G/H$  
is locally isomorphic to a $G$-orbit in the distinguished boundary 
of $\partial_d\Xi$ of $\Xi$. \qed

\subheadline{1.4. The domains $\Xi_H$} 
 
\noindent 
We keep the notation from Subsections 1.2 - 1.3. {}From now on we  
fix an element $X_{\hbox{\fiverm H}}\in \partial_e\Omega$ 
and set $x_{\hbox{\fiverm H}}=\exp(iX_{\hbox{\fiverm H}})$, $z_{\hbox{\fiverm H}}= 
x_{\hbox{\fiverm H}}K_\C$. As the notation suggests we denote by $H<G$ the stabilizer  
of $z_{\hbox{\fiverm H}}\in \partial_d\Xi$ in $G$.  
 
\par For the rest of this paper we will employ the following assumptions: 
$G_\C$ is simply connected and $G/H\simeq G(z_{\hbox{\fiverm H}})$ is an NCC symmetric space (cf.\  
Proposition 1.3.1). Notice that this implies in particular $H=G^\tau$.  
Recall the element $T_0$ from Definition 1.1.1 and notice that we have  
(up to sign) $X_{\hbox{\fiverm H}}={\pi\over 2}T_0$.

\par We define  a domain $\Omega_{\hbox{\fiverm H}}\subeq \Omega$ 
by 
 
$$\Omega_{\hbox{\fiverm H}}=\Int \left(\conv\{{\cal W}(X_{\hbox{\fiverm H}})\}\right)\ .$$ 
Here $\conv\{\cdot\}$ denotes the convex hull of $\{ \cdot\}$ and $\Int(\cdot)$ denotes 
the interior of $(\cdot)$. {}From the definition we immediately obtain that: 
 
\msk 
\item{(1.4.1)} $\Omega_{\hbox{\fiverm H}}$ is open in $\a$. 
\item{(1.4.2)} $0\in \Omega_{\hbox{\fiverm H}}$ (because $ X_{\hbox{\fiverm H}}\neq 0$ 
and ${\cal W}(X_{\hbox{\fiverm H}})$ meets every Weyl chamber). 
\item{(1.4.3)} The set of extremal points of $\oline{\Omega_{\hbox{\fiverm H}}}$ is ${\cal W}(X_{\hbox{\fiverm H}})$. 
\msk 
 
Let us illustrate the geometry for one example. 
 
\Example 1.4.1. Let $G=\Sl(3,\R)$. Then $\a$ is two-dimensional  
and $\Sigma$ is a root system of type $A_2$. We have $\partial_e\Omega={\cal W}(X_1)\amalg {\cal W}(X_2)$ and the  
corresponding isotropy subgroups are given by $H_1=\SO(1,2)$ and $H_2=\SO(2,1)$.  
With $H=H_1$ the geometry of $\Omega$ and $\Omega_H$ is depicted as follows: 
 
\gsk  
\epsfbox{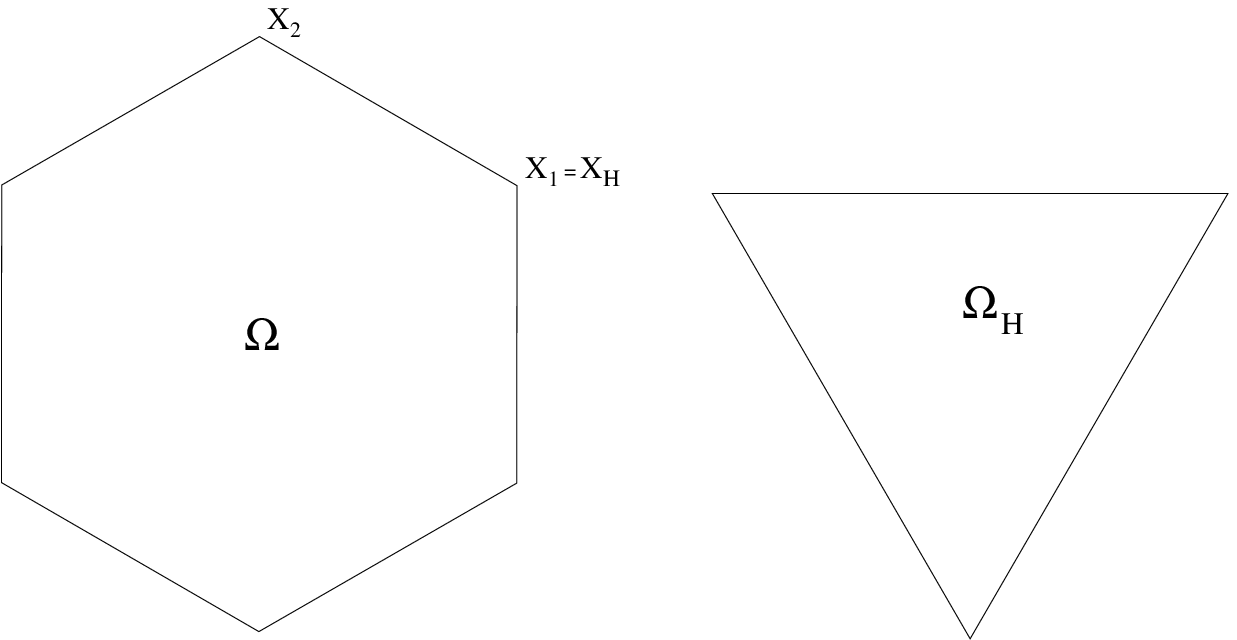}
\gsk\gsk 
\qed 
 
\nin Let us define a domain $\Xi_H\subeq G_\C/ K_\C$ by 
$$\Xi_H=G\exp(i\Omega_{\hbox{\fiverm H}})K_\C/ K_\C\, .$$ 
The domain $\Xi_H$ has the following properties: 
\msk 
\item{(1.4.4)} $\Xi_H$ is $G$-invariant (clear from the definition). 
\item{(1.4.5)} $\Xi_H$ is open in $G_\C/ K_\C$ (follows from (1.4.1) and [AG90]). 
\item{(1.4.6)} $G/K\subeq \Xi_H$ is a totally real submanifold (follows from (1.4.2) and (1.4.5)). 
\item{(1.4.7)} $\Xi_H$ is Stein (follows from the convexity of $\Omega_{\hbox{\fiverm H}}$ and [GK02b]). 
\item{(1.4.8)} $\Xi_H\subeq \Xi$ (because $\Omega_{\hbox{\fiverm H}}\subeq \Omega$). 
\item{(1.4.9)} $\Xi_H=\Xi$ iff $\Sigma$ is of type $C_n$ (cf.\ [KS01b]). 
\msk 
 
\nin  Write $\oline{\Xi_H}$ for the closure of $\Xi_H$ in $G_\C/ K_\C$ and define the 
{\it distinguished boundary} of $\Xi_H$ by 
$$\partial_d\Xi_H=G(z_{\hbox{\fiverm H}})\, .$$ 
Notice that $\partial_d\Xi_H\simeq G/H$ as $G$-spaces. Let us remark further that 
$\partial_d\Xi_H\subeq \partial_d\Xi$.  
 
\par The distinguished boundary $\partial_d\Xi_H$ of $\Xi_H$ can be 
considered as some sort of Shilov boundary of $\Xi_H$. More precisely, 
mimicking the argument in [GK02a, Th.\ 2.3] we obtain that 
 
$$\sup_{z\in \Xi_H} |f(z)|=\sup_{z\in \partial_d\Xi_H|} |f(z)|\leqno (1.4.10)$$ 
for all bounded holomorphic functions $f$ on $\Xi_H$ which continuously extend to $\oline{\Xi_H}$.

\sectionheadline{2. Holomorphic $H$-spherical distributions} 
 
\noindent 
In the section we assume that $G/H$ is NCC symmetric space 
realized as $G(z_{\hbox{\fiverm H}})$ in the distinguished boundary of $\Xi$. 
\par Recall 
that a representation $(\pi ,V)$ of $G$ is called {\it admissible} if  
the multiplicity of each $K$-type is finite and of {\it finite length} if  
the associated Harish-Chandra module of $K$-finite vectors $V_K$ is of finite length. 
Our aim in this section is to 
associate to a non-zero $K$-fixed vector $v_{\hbox{\fiverm K}}$ in an admissible  
representation $(\pi,V)$ of finite length a  certain canonical $H$-spherical 
distribution vector $v_{\hbox{\fiverm H}}\in (V^{-\infty})^H$. For irreducible representations 
$\pi$ the vector $v_{\hbox{\fiverm H}}$ 
is unique in the sense that it allows analytic continuation 
of generalized matrix coefficients on $G/H$ to holomorphic 
functions on $\Xi_H$.  
\par We let $X_{\hbox{\fiverm H}},x_{\hbox{\fiverm H}}$, and $z_{\hbox{\fiverm H}}$ be as in 
the last subsection and recall that $x_{\hbox{\fiverm H}}^{-1}H_\C x_{\hbox{\fiverm H}}=K_\C$. 
 
\msk\nin {\bf General Remark:} In all results of this section  which involve the domain 
$\Xi_H$ one can replace $\Xi_H$ by the bigger domain $\Xi$. This holds in particular  
for the results in Subsection 2.2.

\subheadline{2.1. The definition of the holomorphic $H$-spherical distribution vector}

\noindent 
Before we discuss the general case let us assume for the moment 
that $V$ is irreducible and finite dimensional. Then 
$(\pi, V)$ extends to a holomorphic representation of $G_\C$ 
which we also denote by $(\pi, V)$. 
If $L<G$ is a subgroup of $G$ then we write $V^L$ for the subspace of 
$V$ fixed by $L$. 
Then the mapping 
$$V^K\to V^H, \ \ v_{\hbox{\fiverm K}}\mapsto v_{\hbox{\fiverm H}}\:=\pi(x_{\hbox{\fiverm H}})v_{\hbox{\fiverm K}}$$ 
sets up an isomorphism between the $K$-spherical and $H$-spherical 
vectors of $V$. 
The obvious problem in the general case is, that $\pi(x_{\hbox{\fiverm H}})v$ is 
not necessarily defined as an element in $V$. 
\par 
 
We are now going to develop an appropriate generalization of the mapping 
$v_{\hbox{\fiverm K}}\mapsto v_{\hbox{\fiverm H}}$ for an admissible finite length  
representation $(\pi ,{\cal H})$ of $G$ in a Hilbert space ${\cal H}$. 
Denote by ${\cal H}^\infty$ respectively ${\cal H}^\omega$ the space of smooth respectively 
analytic vectors in ${\cal H}$. Their strong anti-duals, i.e., the space 
of continuous conjugate linear maps into $\C$, are denoted by ${\cal H}^{-\infty}$ 
respectively ${\cal H}^{-\omega}$ and referred 
to as the $G$-modules of {\it distribution 
vectors}, respectively {\it hyperfunction vectors} of $(\pi, {\cal H})$. 
Notice the chain of continuous inclusions ${\cal H}^{\omega}\into {\cal H}^{\infty}\into {\cal H}\into 
{\cal H}^{-\infty}\into {\cal H}^{-\omega}$. Here the  inclusion 
${\cal H}\into {\cal H}^{-\infty}$ is the natural one, $v\mapsto (u\mapsto \la v,u\ra)$. We 
will also use the notation $u\mapsto \la v,u\ra$ for $v\in {\cal H}^{-\infty}$. Define a representation 
$\pi^0$ of $G$ on $\cal H$ by  
$$\eqalign{ 
\la \pi^0(g) v,u\ra &=\la \pi (g^{-1})^*v,u\ra \cr 
&= \la v,\pi (g^{-1})u\ra\cr}\qquad (g\in G; \ u, v\in {\cal H})\ .$$ 
Hence the natural representation $\pi^{-\infty}$ of $G$ on ${\cal H}^{-\infty}$ is an extension of the representation 
$\pi^0$ on ${\cal H}$ to ${\cal H}^{-\infty}$. The representation $\pi^0$ is 
called the {\it conjugate dual representation} of $\pi$. 
Notice that $(\pi^0, {\cal H})$ is admissible and of finite length if and 
only if the same holds for $(\pi ,{\cal H})$. Note also that $\pi =\pi^0$ if and only if 
$\pi $ is unitary. Notice that if $u,v\in {\cal H}$ then $\la v ,u\ra =\overline{\la u, v\ra}$. Accordingly, 
if $u\in {\cal H}^{\infty}$ and $v\in {\cal H}^{-\infty}$, then we write 
$$\la u,v\ra :=\overline{\la v,u\ra}\, .$$

\par 
Denote by ${\cal H}_K$ the $(\g,K)$-module of $K$-finite vectors in ${\cal H}$. Note that by 
our assumption that ${\cal H}$ is admissible and of finite length it follows that ${\cal H}_K\subseteq {\cal H}^\omega$. 
But usually  we cannot find $H$-fixed vectors in ${\cal H}$ but only 
in the larger space of distribution vectors. Another 
complication arises as the space of $({\cal H}^{-\infty})^H$ of $H$-invariants in 
${\cal H}^{-\infty}$ is finite dimensional but in general not one-dimensional. 
For ``generic'' principal series representations of $G$ one has  
$$\dim ({\cal H}^{-\infty})^H=\left| {\cal W}/ {\cal W}_0\right|$$ 
with ${\cal W}_0=N_{K\cap H}(\a)/ Z_{K\cap H}(\a)$ the {\it little Weyl group}. 
Our correspondence will be that we associate to a 
$K$-fixed vector $v_{\hbox{\fiverm K}}$ a unique $H$-fixed distribution vector $v_{\hbox{\fiverm H}}$. 
If $\pi$ is irreducible, then this 
distribution $v_{\hbox{\fiverm H}}$ is the (up to scalar) unique element of  $({\cal H}^{-\infty})^H$ 
for which generalized matrix coefficients extend to holomorphic functions 
from $G/H$ to $\Xi$. For the proof we will need the following {\it 
Automatic Continuity} Theorem of van den Ban, Brylinski and Delorme  (c.f. [vdBD88, Th. 2.1] and  
[BD92, Th.\ 1] for the version used here). 
\par If $V$ is a complex vector space, then  let us denote by 
$V^{\oline *}$ its algebraic anti-dual. Set 
$({\cal H}^{-\infty})^{\h}=\{v\in {\cal H}^{-\infty}\: (\forall X\in \h)\, d\pi^{-\infty}(X)v=0\}$.

\Theorem 2.1.1. {\bf (Automatic Continuity)} 
Let $(\pi, {\cal H})$ be an admissible representation of $G$ with finite length. Then 
$$({\cal H}^{-\infty})^{H_0}\simeq ({\cal H}_K^{\oline *})^\h\ , \leqno (2.1.1) $$ 
meaning every $\h$-fixed anti-linear functional on ${\cal H}_K$ admits 
a unique extension to a continuous and $H_0$-fixed anti-linear functional 
on ${\cal H}^\infty$. In particular we have that 
$$({\cal H}^{-\infty})^{H_0}\simeq ({\cal H}^{-\omega})^{H_0}\ .\ \leqno(2.1.2)$$\qed 
 
\par Recall the 
$G-K_\C$ double coset domain $\tilde \Xi=G\exp(i\Omega)K_\C$ 
in $G_\C$ and the complex crown  $\Xi=\tilde\Xi/ K_\C$. Our methods use holomorphic extensions of 
representations. We recall therefore some results from [KS01a] \S 4 and in particular Theorem 3.1: 
 
\Theorem 2.1.2.  Let $(\pi, {\cal H})$ be an admissible Hilbert representation of $G$ 
with finite length. Then the following assertions hold: 
 
\item{(i)} For every $v\in {\cal H}_K$ the orbit 
mapping $G\to {\cal H}, \ g\mapsto \pi(g)v$ extends to a $G$-equivariant holomorphic 
mapping $\tilde \Xi\to {\cal H}$. 
 
\item{(ii)} Let $v,w\in {\cal H}_K$. Then the restricted matrix coefficient 
$A\to \C, \ a\mapsto \la \pi (a)v, w\ra $ extends to a holomorphic mapping to 
the abelian tube domain $T (2\Omega)=A\exp(2i\Omega)$.\qed 
 
We will from now on assume that $\pi\res_K$ is unitary. This is no restriction 
in view of Weyl's unitarity trick. 
 
\par 
As before  we identify $G/H$ with the subset $G(z_{\hbox{\fiverm H}})$ of $\overline{\Xi_H}$. 
For $0\leq t<1$ we set $a_t=\exp(itX_{\hbox{\fiverm H}})\in \exp(i\Omega)$ and notice that 
$\lim_{t\to 1} a_t =x_{\hbox{\fiverm H}}$. 
 
\Theorem 2.1.3. Let $(\pi, {\cal H})$  be a $K$-spherical admissible  representation 
of $G$ with finite length. Let $v_K\in {\cal H}^K$ be non-zero. Then the anti-linear functional 
 
$$v_{\hbox{\fiverm H}}^\omega\: {\cal H}^\omega \to \C , \ \ v\mapsto \lim_{t\nearrow 1} \la \pi^0(a_t)v_{\hbox{\fiverm K}}, v\ra$$ 
is well defined, non-zero and admits a non-trivial extension to an  
$H$-fixed distribution vector of $(\pi, {\cal H})$. \qed 
 
\ssk \nin {\bf Note:} If $(\pi, {\cal H})$ is unitary, then $\pi^0(a_t)=\pi(a_t)$ and so 
$v_{\hbox{\fiverm H}}^\omega(v)=\lim_{t\nearrow 1} \la \pi(a_t)v_{\hbox{\fiverm K}}, v\ra$ for all $v\in {\cal H}^\omega$.

\Proof.  We first show that  $v_{\hbox{\fiverm H}}^\omega$ is well defined.  
Let $v\in {\cal H}^\omega$. As $v$ is  
analytic,  we find an $0<\epsilon <1$ such that  
$\pi (a_\eps)v$ is defined. Notice that $\pi^0 (a_{1-\eps})v_{\hbox{\fiverm K}}$ is defined 
by Theorem 2.1.2 (i).  
Then we have for all $\eps<t<1$ that  
$$\la \pi^0(a_t)v_{\hbox{\fiverm K}}, v\ra= \la \pi^0(a_{t-\eps})v_{\hbox{\fiverm K}}, \pi(a_\eps)v\ra$$ 
and hence 
$$\lim_{t\nearrow 1}\la \pi^0(a_t)v_{\hbox{\fiverm K}}, v\ra=\la \pi^0(a_{1-\eps})v_{\hbox{\fiverm K}},  
\pi(a_\eps)v\ra\ .$$ 
Thus  $v_{\hbox{\fiverm H}}^\omega$ is defined. 
 
\par  Next we show that $v_{\hbox{\fiverm H}}^\omega$ is fixed by $H$.  
As $G/H$ is NCC, it follows that $H=H_0 Z_{H\cap K}(\a )$ [H\'O96, p.79]. Further it is clear from the definition  
that $v_{\hbox{\fiverm H}}^\omega$it is fixed by $Z_{H\cap K}(\a )$. Hence it is  
enough to prove that $v_{\hbox{\fiverm H}}^\omega$ 
is $H_0$-fixed, i.e. annihilated by $\h$. For that let $Y\in \h$ and 
$v\in {\cal H}^\omega$. Then 
 
$$\eqalign{v_{\hbox{\fiverm H}}^\omega(d\pi(Y)v)&=\lim_{t\nearrow 1}\la \pi^0(a_t)v_{\hbox{\fiverm K}}, d\pi(Y)v\ra= 
-\lim_{t\nearrow 1}\la d\pi^0(Y)\pi^0(a_t)v_{\hbox{\fiverm K}}, v\ra\cr 
&= - \lim_{t\nearrow 1}\la\pi^0(a_t)d\pi^0(\Ad(a_t)^{-1} Y)v_{\hbox{\fiverm K}}, v\ra\cr 
&=  -\lim_{t\nearrow 1}\la\pi^0(a_{t-\eps})d\pi^0(\Ad(a_t)^{-1} Y)v_{\hbox{\fiverm K}}, \pi(a_\eps)v\ra\ .\cr}$$ 
{}From Theorem 2.1.2 (i) we now obtain that 
$$\pi^0(a_{t-\eps})d\pi^0(\Ad(a_t)^{-1} Y)v_{\hbox{\fiverm K}}\to 
\pi^0(a_{1-\eps})d\pi^0(\Ad(x_{\hbox{\fiverm H}})^{-1} Y)v_{\hbox{\fiverm K}}\ .$$ 
But $\Ad (x_{\hbox{\fiverm H}})^{-1}Y\in \k_\C$ and hence 
$d\pi^0(\Ad (x_{\hbox{\fiverm H}})^{-1}Y)v_{\hbox{\fiverm K}}=0$. 
We thus get: 
$$\eqalign{v_{\hbox{\fiverm H}}^\omega(d\pi(Y)v) 
&=-\lim_{t\nearrow 1}\la\pi^0(a_{t-\eps})d\pi^0(\Ad(a_t)^{-1} Y)v_{\hbox{\fiverm K}}, \pi(a_\eps)v\ra\cr 
&= -\la\pi^0(a_{1-\eps})d\pi^0(\Ad(x_{\hbox{\fiverm H}})^{-1} Y)v_{\hbox{\fiverm K}}, \pi(a_\eps)v\ra\cr 
&= 0\, .\cr}$$

\par Finally, let us show that $v_{\hbox{\fiverm H}}^\omega\neq 0$. Suppose the contrary, i.e.  
$v_{\hbox{\fiverm H}}^\omega=0$. Then it follows  that  
$$(\forall u\in {\cal H}_K) \qquad \la v_{\hbox{\fiverm H}}^\omega, u\ra =0\ . \leqno (2.1.3)$$ 
Now consider the function  
$$f\: \R +i]-2,2[\to \C , \ \ z\mapsto \la \pi^0(\exp(zX_{\hbox{\fiverm H}}))v_{\hbox{\fiverm K}},  
v_{\hbox{\fiverm K}}\ra\ .$$ 
According to Theorem 2.1.2 (ii) the function $f$ is well defined and holomorphic.  
It is clear that $f\not\equiv 0$ as $f(0)=\la  v_{\hbox{\fiverm K}},v_{\hbox{\fiverm K}}\ra>0$. 
But (2.1.3) implies that $f^{(n)}(i)=0$ for all $n\in \N_0$; a contradiction  
to $f\not\equiv 0$. \qed  
 
In the sequel we write $v_{\hbox{\fiverm H}}$ for the $H$-fixed distribution vector 
obtained from $v_{\hbox{\fiverm H}}^\omega$. We will call 
$v_{\hbox{\fiverm H}}$ the {\it holomorphic $H$-spherical distribution vector} of $(\pi, {\cal H})$ 
corresponding to $v_{\hbox{\fiverm K}}$.

\Remark 2.1.4. By Theorem 2.1.3 we have 
$$v_{\hbox{\fiverm H}}^\omega=\hbox{w}-\lim_{t\nearrow 1} \pi^0(a_t)v_{\hbox{\fiverm K}}, $$ 
i.e. $v_{\hbox{\fiverm H}}^\omega$ is the weak limit of $\pi^0(a_t)v_{\hbox{\fiverm K}}$ for $t\to 1$ in the 
locally convex space ${\cal H}^{-\omega}$.  
It is possible to strengthen this convergence: Let $v\in {\cal H}^\omega$ and $C\subeq G$  
a compact subset. Then for all $\eps>0$ there exists $0<s<1$ such that for all 
$0<s<t<1$: 
 
$$\sup_{g\in C} |\la v_{\hbox{\fiverm H}}^\omega, \pi(g)v\ra -\la \pi^0(a_t)v_{\hbox{\fiverm K}}, \pi(g)v\ra|<\eps\ 
 . \leqno(2.1.4)$$ 
In fact this follows from a simple modification of the first part of the  
proof of Theorem 2.1.3: we only have to observe that for $v\in {\cal H}^\omega$ 
there exists a $0<\delta<1$ such that $\pi(a_\delta)\pi(g)v$ exists for all  
$g\in C$.

\par Supported by calculations in the rank one case we conjecture 
that one actually has 
$ \pi(a_t)v_{\hbox{\fiverm K}}\to v_{\hbox{\fiverm H}}$ in ${\cal H}^{-\infty}$ weakly (and hence strongly by the 
Banach-Steinhaus Theorem which applies as ${\cal H}^\infty$ is a Fr\'echet space ). \qed 
 
Let us illustrate the situation by the discussion of one example. 
 
\Example 2.1.5. Here we will determine an explicit analytic description 
of $v_{\hbox{\fiverm H}}$ for unitary principal series of $G=\Sl(2,\R)$. 
Let $(\pi_\lambda, {\cal H}_\lambda)$ denote a unitary spherical 
principal series of $G$ with parameter $\lambda\in i\a^*$. 
Then $\pi_\lambda^0=\pi_\lambda$ as $\pi_\lambda$ is unitary.  In the sequel 
we will identify $\a_\C^*$ with $\C$ in such a way that $\rho\in \a^*$ corresponds to $1$. 
With our choice of $\a$ to be 
$$\a=\left\{ \pmatrix{s & 0\cr 0& -s\cr}\: s\in \R\right\}$$ 
this identification is given by 
$$\lambda \mapsto \lambda \pmatrix{1 & 0 \cr 0 & -1\cr}\, .$$ 
 
\par We will use the noncompact realization of ${\cal H}_\lambda =L^2(\R)$ of 
$\pi_\lambda$. Then for $g=\pmatrix {a & b\cr c& d\cr}\in G$ the operator 
$\pi_\lambda(g)$ is given by 
 
$$(\pi_\lambda(g)f)(x)=|bx+d|^{-1-\overline{\lambda}} f\left({ax+c\over bx +d}\right) 
\qquad (f\in L^2(\R), x\in \R)\ .$$ 
A normalized $K$-spherical vector is then given by 
 
$$v_{\hbox{\fiverm K}}(x)={1\over \sqrt\pi} (1+x^2)^{-{1\over 2}(1+\overline{\lambda})}\ .$$ 
\par Notice that $H=\SO(1,1)$ and  (up to sign) we have  
$X_{\hbox{\fiverm H}}=\pmatrix {\pi\over 4 & 0\cr 0& -{\pi\over 4}\cr}$.  
Thus for $0\leq t<1$ the element $a_t\in 
\exp(i\Omega)$ is given by  
 
$$a_t=\pmatrix{ e^{i{\pi\over 4}t} & 0 \cr 0 & e^{-i{\pi\over 4}t} \cr}\  .$$ 
Then we have 
$$v_{\hbox{\fiverm H}}^\omega=\hbox{w}-\lim_{t\nearrow 1} \pi_\lambda(a_t)v_{\hbox{\fiverm K}}$$ 
or 
$$v_{\hbox{\fiverm H}}^\omega(x)=\hbox{w}- \lim_{t\nearrow 1}{e^{i{\pi\over 4}t 
(1+\overline{\lambda})}\over \sqrt\pi} (1+e^{i\pi t}x^2)^{-{1\over 2}(1+\overline{\lambda})}\ .$$ 
A simple calculation then shows that $v_{\hbox{\fiverm H}}^\omega$ and $v_{\hbox{\fiverm H}}$ 
are given by the locally integrable function 
 
$$v_{\hbox{\fiverm H}} (x)=\cases {{e^{i{\pi\over 4}(1+\overline{\lambda})}\over \sqrt\pi} 
(1- x^2)^{-{1\over 2}(1+\overline{\lambda})} & for $|x|<1$, \cr 
0 & for $|x|=1$, \cr 
{e^{-i{\pi\over 4}(1+\overline{\lambda})} \over \sqrt \pi} 
(x^2-1)^{-{1\over 2}(1+\overline{\lambda})} & for $|x|>1$\ .\cr} $$ 
A basis of $({\cal H}_\lambda^{-\infty})^H$ is 
given by $v_{\hbox{\fiverm H,1}}, v_{\hbox{\fiverm H, 2}}$ where 
$$v_{\hbox{\fiverm H,1}}(x)=\cases { {1\over \sqrt\pi} (1- x^2)^{-{1\over 2}(1+\overline{\lambda})} & for $|x|<1$, \cr 
0 & for $|x|\geq 1$. \cr}$$ 
and 
$$v_{\hbox{\fiverm H,2}}(x)=\cases {  {1\over \sqrt\pi}(x^2-1)^{-{1\over 2}(1+\overline{\lambda})} & for $|x|>1$, \cr 
0 & for $|x|\leq 1$. \cr}$$ 
These two basis vectors are chosen such  that they have support in the 
open $H$-orbits on $\P^1(\R)$, namely $]-1, 1[$ and $\P^1(\R)\bs [-1, 1]$. 
Notice that $v_{\hbox{\fiverm H}}$ is a non-trivial linear combination of $v_{\hbox{\fiverm H, 1}}$ and  
$v_{\hbox{\fiverm H,2}}$ and that 
$v_{\hbox{\fiverm H}}$ has full support on $\R$. 
Another interesting feature of this example is that we  have here 
$$v_{\hbox{\fiverm H}} =\hbox{w}-\lim_{t\nearrow 1}  
\pi_\lambda(a_t)v_{\hbox{\fiverm K}} \qquad \hbox{in ${\cal H}_\lambda^{-\infty}$}$$ 
and hence also strongly by the Banach-Steinhaus Theorem (compare with the conjecture  
stated at the end of Remark 2.1.4) . \qed

\subheadline{2.2. Holomorphic extension of matrix coefficienst on $G/H$ to $\Xi_H$} 
 
\noindent 
In this subsection we clarify the role of the holomorphic distribution vector 
$v_{\hbox{\fiverm H}}$ in view of holomorphic extensions of matrix coefficients from 
$G/H$ to $\Xi$.  If $U$ is a complex manifold, then we denote by 
${\cal O}(U)$ the space of holomorphic functions $f:U\to \C$. 
 
\Definition 2.2.1. Let $f$ be a continuous function on $G/H$. Then we say 
that $f$ has a {\it holomorphic extension} to $\Xi_H$ if there exists an 
$\tilde f\in {\cal O}(\Xi_H)$ such that for all compact subsets 
$C\subeq G$  one has 
$$\lim_{t\nearrow 1} \sup_{g\in C} |f(gH) -\tilde f(ga_tK_\C)|=0\ .\leqno(2.2.1)$$ 
\qed 
 
Notice that (2.2.1) implies that 
$$f(gH)=\lim_{t\nearrow 1} \tilde f(ga_tK_\C)$$ 
for all $g\in G$. 
Furthermore the holomorphic extensions are unique by  the following lemma:

\Lemma 2.2.2. {\rm\bf (Identity Theorem for holomorphic extensions)} Let $f\in C(G/H)$ 
and assume that $f$ has 
a holomorphic extension $\tilde f\in {\cal O}(\Xi_H)$. Then $f\equiv 0$ implies 
$\tilde f\equiv 0$. In particular, the holomorphic extension $f\in C(G/H)$ 
is unique if it exists. 
 
\Proof. This is easily reduced to the one-dimensional case as 
follows. Define an abelian tube domain 
$T=\exp(\R X_{\hbox{\fiverm H}} + ]-1,1[ iX_{\hbox{\fiverm H}})$ and set $\partial_s 
T =\exp(\R X_{\hbox{\fiverm H}} + iX_{\hbox{\fiverm H}})$. We realize 
$T\subeq \Xi_H$ and $\partial_sT\subeq G/H$ through the $T$-orbit , respectively $\partial_s T$-orbit, 
through $K_\C\in \Xi_H$, respectively $z_{\hbox{\fiverm H}}\in G/H$. 
Let $f\in C(G/H)$ and assume that $f$ has a holomorphic extension $\tilde f$. 
Let $\phi=f\res_{\partial_s T}$. Then $\phi$ has a holomorphic extension to 
$T$ given by $\tilde \phi=\tilde f\res_T$. By the well known one-dimensional 
situation we have $\tilde\phi\equiv 0$ if $\phi\equiv 0$. Thus if $f\equiv 0$, we obtain 
$\tilde f\res_T\equiv 0$. Replacing $f$ by $f_g$ where $f_g(xH)=f(gxH)$ for $g\in G$, then 
the above discussion implies that 
$\tilde f\res_{GT}\equiv 0$ and hence $\tilde f\equiv 0$ as $GT$ contains the totally 
real submanifold $G/K$ of $\Xi_H$. \qed 
 
We assume that the complex conjugation $\g_\C\to\g_\C, \ X\mapsto \oline X$ with respect to 
the real form $\g$ of $\g_\C$ lifts to a conjugation $g\mapsto \oline g$ of $G_\C$. 
Notice that this is always satisfied if $G_\C$ is the universal complexification of $G$. 
Notice also that $\overline{x}\in \tilde{\Xi}$ for all $x\in\tilde{\Xi}$. 
Let $u,v\in {\cal H}_K$. Then the functions 
$$\tilde{\Xi}\ni x\mapsto\la \pi^0(x)u,v\ra , \la u,\pi(\overline{x}^{-1})v \ra\in \C $$ 
are well defined and holomorphic by Theorem 2.1.2. Both of them agree if $x\in G$ and 
hence they agree on all of $\tilde{\Xi}$. 
 
\Proposition 2.2.3. Let $(\pi, {\cal H})$ be a admissible $K$-spherical representation of $G$ 
with finite length. Let $v\in {\cal H}^\omega$. Then the following assertions hold: 
 
\item{(i)} The 
matrix coefficient 
$$f_{v, v_{\hbox{\fiverm H}}}\: G/H\to \C, \ \ gH\mapsto \la \pi (g^{-1})v, v_{\hbox{\fiverm H}}\ra$$ 
admits a holomorphic extension $\tilde f_{v, v_{\hbox{\fiverm H}}}$. 
Moreover, we have 
$$\tilde f_{v, v_{\hbox{\fiverm H}}}(xK_\C)=\la v, \pi^0(\overline x)v_{\hbox{\fiverm K}}\ra $$ 
for all $xK_\C\in \Xi_H$. 
\item{(ii)} The 
matrix coefficient 
$$g_{v_{\hbox{\fiverm H}},v }\: G/H\to \C, \ \ gH\mapsto \la \pi^0(g)v_{\hbox{\fiverm H}},v\ra$$ 
admits a holomorphic extension $\tilde g_{v_{\hbox{\fiverm H}},v}$. 
Moreover, we have 
$$\tilde g_{v_{\hbox{\fiverm H}},v}(xK_\C)=\la \pi^0(x )v_{\hbox{\fiverm K}},v\ra \quad\hbox{\rm and} 
\quad \tilde{g}_{v_{\hbox{\fiverm H}},v}(xK_\C) 
=\overline{\tilde f_{v, v_{\hbox{\fiverm H}}}(\overline{x}K_\C)}$$ 
for all $xK_\C\in \Xi_H$.
 
\Proof. We will only show (i) as the proof for (ii) is the same. By Theorem 2.1.2 (i) it follows that 
$$\tilde{f}(xK_\C)=\la v,\pi^0(\bar{x})v_{\hbox{\fiverm K}}\ra$$ 
exists and is holomorphic on $\Xi_H$. Let $g\in G$ and $0<t<1$. Then 
$$\la \pi (g^{-1})v,\pi^0(a_t)v_{\hbox{\fiverm K}}\ra =\la v,\pi^0(ga_t)v_{\hbox{\fiverm K}}\ra 
=\tilde{f}(ga_tK_\C)\, .$$ 
Taking the limit at $t\to 1$ and using the remarks just before (2.1.4) it follows that 
$$\lim_{t\nearrow 1} \tilde{f}(ga_tK_\C)=\la \pi (g^{-1})v,v_{\hbox{\fiverm H}}\ra$$ 
and the convergence is uniform on compact subsets in $G$. 
\qed 
 
Denote by $({\cal H}^{-\infty})^H_{\rm hol}\subset ({\cal H}^{-\infty})^H$ 
the space of $H$-invariant distribution vectors $\eta$ such 
that the function 
$$G/H\ni x \mapsto g_{\eta ,v}(x):=\la \pi^0(x)\eta ,v\ra$$ 
has a holomorphic extension to $\Xi_H$  
for all $v\in {\cal H}^\omega$. Notice that for $g\in G$ and $x\in \Xi_H$ we have 
$$\tilde g_{\eta ,\pi (g)v}(x)=\tilde g_{\eta ,v}(g^{-1}x)\, .\leqno(2.2.2)$$ 
Our next task is to prove a converse of Proposition 2.2.3, namely that 
the map 
$${\cal H}^K\ni v_{\hbox{\fiverm K}}\mapsto v_{\hbox{\fiverm H}}\in ({\cal H}^{-\infty})^H_{\rm hol}$$ 
is an isomorphism. In particular only the $H$-invariant distribution 
vectors constructed in Theorem 2.1.3 have holomorphic extension. Thus if 
$\dim  {\cal H}^K=1$, as in the case of the principal series 
representations of $G$, the space $({\cal H}^{-\infty})^H_{\rm hol}$ is 
also one-dimensional, i.e., there is (up to scalar) a unique $H$-spherical distribution vector which allows 
holomorphic extension of the  smooth matrix coefficients. 
\par 
As before we consider the pairing 
$$ ({\cal H}^{-\infty})^H\times {\cal H}^\infty \to C^\infty(G/H), \ \ (\eta,v)\mapsto 
g_{\eta,v}; \ g_{\eta,v}(gH)=\la \pi^0(g)\eta,v\ra\ .$$

\Theorem 2.2.4. Let $(\pi , {\cal H})$ be an admissible representation of $G$ with finite length. Then the map 
$${\cal H}^K\ni v_{\hbox{\fiverm K}}\mapsto v_{\hbox{\fiverm H}}\in ({\cal H}^{-\infty})^H_{\rm hol}$$ 
is a linear isomorphism. 
 
\Proof. It follows from Theorem 2.1.3 that ${\cal H}^K\ni v_{\hbox{\fiverm K}}\mapsto 
v_{\hbox{\fiverm H}}\in({\cal H}^{-\infty})^H_{\rm hol}$ is well defined and injective. 
It is also clear that the map is linear. It remains to show that the map is onto.  
For that let $\eta\in ({\cal H}^{-\infty})^H_{\rm hol}$. Define a conjugate linear map 
$\tilde\eta : {\cal H}_K\to \C$ by 
$$\tilde\eta(u)= \tilde g_{\eta ,u}(K_\C)\, .$$ 
Then it follows from (2.2.2) that $\tilde\eta$ is $K$-invariant. Thus we find a unique  
$v_{\hbox{\fiverm K}}\in {\cal H}_K$ such that  
$\tilde\eta(u)=\la v_{\hbox{\fiverm K}}, u\ra$  for all $u\in {\cal H}_K$.  
In particular,  it follows that  
$$(\forall u\in {\cal H}_K)\qquad  
\tilde g_{\eta ,u}(K_\C)= \tilde g_{v_{\hbox{\fiverm H}} ,u}(K_\C)\ .\leqno(2.2.3)$$ 
\par Fix $w\in {\cal H}_K$. We claim that $\tilde g_{\eta ,w}=\tilde g_{v_{\hbox{\fiverm H}} ,w}$.  
In fact, it follows from (2.2.2) and (2.2.3) that all derivatives  
of  $\tilde g_{\eta ,w}$ and $\tilde g_{v_{\hbox{\fiverm H}} ,w}$ coincide at $K_\C \in\Xi_H$.  
Thus $\tilde g_{\eta ,w}=\tilde g_{v_{\hbox{\fiverm H}} ,w}$ by Taylor's Theorem.  
 
\par {} It follows from our claim that  $g_{\eta ,w}=g_{v_{\hbox{\fiverm H}} ,w}$ for all 
$w\in {\cal H}_K$. In particular, we obtain that   
$$(\forall w\in {\cal H}_K)\qquad  
\la \eta, w\ra =g_{\eta,w}(H)=g_{v_{\hbox{\fiverm H}}, w}(H)=\la v_{\hbox{\fiverm H}}, w\ra\ , $$ 
and so $\eta=v_{\hbox{\fiverm H}}$, concluding the proof of the theorem.\qed

\Corollary 2.2.5. {\rm\bf (Multiplicity one)}  Let $(\pi, {\cal H})$ be an irreducible $K$-spherical  
Hilbert representation of $G$. Then $\dim ({\cal H}^{-\infty})_{\rm hol}^H =1$.  
 
\Proof. According to [H84, Ch.\ IV, Th. 4.5(iii)] we have that $\dim {\cal H}^K=1$.  
Thus the assertion follows from Theorem 2.2.4.\qed

\subheadline{2.3. Distributional characters and boundary values} 
 
\noindent 
Denote by $dg$ and $dh$  Haar measures on $G$ and $H$. Notice 
that both $G$ and $H$ are unimodular and hence a left Haar measure is 
also a right Haar measure. Denote by $dgH$ a invariant measure on $G/H$, which 
we will normalize in a moment. 
Recall that the mapping  
$$C_c^\infty (G)\to C_c^\infty (G/H), \ \ f\mapsto f^H; \ f^H(xH)=\int_H f(xh)\ dh$$ 
is continuous and onto. We will normalize the measure $dgH$ in such a way that 
$$\int_G f(g) \ dg =\int_{G/H} f^H(gH)\ dgH$$ 
holds for all $f\in C_c(G)$.

\par In this section $(\pi, {\cal H})$ will denote a {\it unitary}  
admissible representation representation of $G$ with finite length.  
Further we will assume that ${\cal H}^K\not=\{0\}$. Notice that  
$\pi$ unitary implies $\pi=\pi^0$.  
\par For $f\in C_c^\infty(G)$ let us recall  
the mollifying property: $\pi(f){\cal H}^{-\infty}\subeq 
{\cal H}^\infty$. Thus the mapping  
$$\Theta_{\pi,v_{\hbox{\fiverm H}}} \: C_c^\infty(G)\to \C, \ \ f\mapsto 
\la \pi(f)v_{\hbox{\fiverm H}}, v_{\hbox{\fiverm H}}\ra\, .$$ 
is well defined. It is known that $\Theta_\pi=\Theta_{\pi,v_{\hbox{\fiverm H}}}$ is a $H$-bi-invariant  
positive definite distribution on $G$. The $H$-bi-invariance implies  
that $\Theta_\pi(f)$ does only depend on $f^H$. We can therefore define a $H$-invariant 
distribution on $G/H$, also denoted by $\Theta_\pi$,  by 
$$\Theta_\pi (f^H)=\Theta_\pi(f)\, .$$ 
 
\par On the other hand we notice that $\pi (x)v_{\hbox{\fiverm K}}\in {\cal H}^\omega$ 
for all $x\in G\exp (i\Omega_H)K_\C$ and hence 
$x\mapsto \la v_{\hbox{\fiverm H}},\pi (x)v_{\hbox{\fiverm K}}\ra$ 
descends to a well defined and anti-holomorphic function on $\Xi_H$.  
Here $v_{\hbox{\fiverm K}}$ and $v_{\hbox{\fiverm H}}$ correspond to each other 
according to Theorem 2.2.4. We can therefore define the holomorphic function 
$\theta_\pi=\theta_{\pi,v_H}:\Xi_H\to \C$ by 
$$\theta_\pi(xK_\C)= \overline{\la v_{\hbox{\fiverm H}}, \pi(x)v_{\hbox{\fiverm K}}\ra }=\la \pi (x)v_{\hbox{\fiverm K}} , 
v_{\hbox{\fiverm H}}\ra \, .$$ 
Then $\theta_\pi$ is left $H$-invariant.

Our next aim is to show that $\Theta_\pi$ is given by the limit 
operation and convolution: 
 
$$\Theta_\pi(f)=\lim_{t\nearrow 1}\int_{G/H} 
f(gH)\ \overline{\theta_\pi(g^{-1}a_t)} \ dgH$$ 
for suitable regular functions $f$. 
 
\par 
As we have only established the convergence $\pi(a_t)v_{\hbox{\fiverm K}}\to v_{\hbox{\fiverm H}}$ in 
${\cal H}^{-\omega}$ and not in ${\cal H}^{-\infty}$, we cannot work 
with test-functions but must use an appropriate  space of 
analytic vectors.  
For that let us write  
$L^1(G)^{\omega, \omega}$ for the space  
of analytic vectors for the left-right regular representation of $G\times G$ on $L^1(G)$. 
Notice that $L^1(G)^{\omega, \omega}$ is an algebra under convolution  
which is  invariant under the natural involution $f\mapsto f^*$ with $f^*(x)=\oline {f(x^{-1})}$. 
Its importance lies in the fact that  $\pi(f){\cal H}^{-\omega}\subeq {\cal H}^\omega$ holds for all  
$f\in  L^1(G)^{\omega,\omega}$ (cf.\ Proposition A.4.1 in the appendix).  
 
\par Write $L^1(G/H)^\omega$ for the space of analytic vectors 
for the left regular representation of $G$ on $L^1(G/H)$. 
According to Proposition A.3.2 below,  the averaging map $f\mapsto f^H$ maps 
$L^1(G)^{\omega,\omega}$ into $L^1(G/H)^\omega$. 
Finally let us define the space: 
$${\cal A}^1(G/H)=\{f^H\in L^1(G/H)^\omega\mid f\in L^1(G)^{\omega,\omega}\}\ .  $$ 
 
\par {}From our discussion above we conclude that the mapping 
$$\Theta_\pi^\omega\: L^1(G)^{\omega,\omega} \to \C, \ \ f\mapsto \la \pi(f)v_{\hbox{\fiverm H}}, v_{\hbox{\fiverm H}}\ra$$ 
is well defined and $H$ bi-invariant. In particular  $\Theta_\pi^\omega(f)$ depends only on $f^H$ and therefore 
factors to ${\cal A}^1(G/H)$. We denote the corresponding map again by $\Theta_\pi^\omega$.

\Theorem 2.3.1. Let $(\pi,{\cal H})$ be an unitary  admissible  representation 
of $G$ of finite length with ${\cal H}^K\not=\{0\}$. Let $v_{\hbox{\fiverm K}}\in {\cal H}^K$. Then  
 
$$(\forall f\in {\cal A}^1(G/H))\qquad \Theta_\pi^\omega(f)=\lim_{t\nearrow 1}\int_{G/H} 
f(gH)\ \overline{\theta_\pi(g^{-1}a_t)}\ dgH\ .$$ 
 
\Proof. Let $F\in L^1(G)^{\omega,\omega}$ be such that $f=F^H$. 
As $\pi(F)v_{\hbox{\fiverm H}}\in {\cal H}^\omega$ (cf.\ Proposition A.4.1), it follows from 
Theorem 2.1.3 that 
$$\Theta_\pi^\omega(F)=\lim_{t\nearrow 1}\la \pi(F)v_{\hbox{\fiverm H}}, \pi(a_t)v_{\hbox{\fiverm K}}\ra\ .$$ 
Thus: 
 
$$\eqalign{\Theta_\pi^\omega(F) 
&=\lim_{t\nearrow 1}\la v_{\hbox{\fiverm H}}, \int_G \oline{F(g^{-1})}\pi(g)\pi(a_t)v_{\hbox{\fiverm K}}\ dg\ra\cr 
&=\lim_{t\nearrow 1}\int_G F(g^{-1})\  \la v_{\hbox{\fiverm H}},\pi(g)\pi(a_t)v_{\hbox{\fiverm K}}\ra \ dg\cr 
&=\lim_{t\nearrow 1}\int_G F(g)\  \la v_{\hbox{\fiverm H}},\pi(g^{-1})\pi(a_t)v_{\hbox{\fiverm K}}\ra \ dg\cr 
&=\lim_{t\nearrow 1}\int_G  F(g) \ \overline{\theta_\pi(g^{-1}a_t)}\ dg\ .\cr}\leqno (2.3.1)$$ 
Fix $0\leq t<1$. We claim that  $g\mapsto \theta_\pi(g^{-1}a_t)$ is a bounded  
function on $G$. Indeed, for $g\in G$  we have  
$$\eqalign{|\theta_\pi(g^{-1}a_t)|&= 
|\la \pi(g^{-1}a_t)v_{\hbox{\fiverm K}},v_{\hbox{\fiverm K}}\ra|= 
|\la \pi(a_t)v_{\hbox{\fiverm K}},\pi(g)v_{\hbox{\fiverm K}}\ra|\cr  
&\leq \|\pi(a_t)v_{\hbox{\fiverm K}}\|\cdot 
\|\pi(g)v_{\hbox{\fiverm K}}\|\leq \|\pi(a_t)v_{\hbox{\fiverm K}}\|\cdot \| v_{\hbox{\fiverm K}}\|\cr}$$ 
and our claim follows from $\|\pi(a_t)v_{\hbox{\fiverm K}}\|<\infty$.  
Combining (2.3.1) with our claim then yields  
$$\eqalign{\Theta_\pi^\omega(F)&= 
\lim_{t\nearrow 1}\int_G  F(g) \ \overline{\theta_\pi(g^{-1}a_t)}\ dg\cr 
& =\lim_{t\nearrow 1}\int_{G/H} f(gH)\  \overline{\theta_\pi(g^{-1}a_t)}\ dgH\ ,\cr}$$ 
completing the proof of the theorem. \qed 
 
We finish this subsection by the following simple remark. 
 
\Lemma 2.3.2. Suppose that $\pi$ is irreducible. Then $\Theta_\pi$ is an eigendistribution 
of the algebra $\D(G/H)$ of invariant differential operators on $G/H$. 
 
\Proof. Recall the surjective homomorphisms ${\cal U}(\g_\C)^{\frak h} 
\to \D(G/H)$ 
and ${\cal U}(\g_\C)^{\frak k}\to \D(G/K)$. We have $x_{\hbox{\fiverm H}}^{-1}H_\C x_{\hbox{\fiverm H}} 
=K_\C$. Hence, ${\rm Ad}(x_{\hbox{\fiverm H}}^{-1})$ defines an isomorphism 
$${\rm Ad}(x_{\hbox{\fiverm H}}^{-1})\:  
{\cal U}(\g_\C)^{\frak h}\to {\cal U}(\g_\C)^{\frak k}\, .$$ 
\par In order to prove the lemma it is sufficient to show  
that $v_{\hbox{\fiverm H}}$ is an eigenvector for each  
$d\pi^0(u)$, $u\in  {\cal U}(\g_\C)^{\frak h}$.  
\par Notice that for each $\tilde u\in  {\cal U}(\g_\C)^{\frak k}$ 
there exists a constant $c(\tilde u)$ such that  
$d\pi^0(\tilde u)v_{\hbox{\fiverm K}}  
= c(\tilde u) v_{\hbox{\fiverm H}}$.  
For $u\in{\cal U}(\g_\C)^{\frak h}$ we now obtain that 
 
$$\eqalign{d\pi^0(u)v_{\hbox{\fiverm H}}&= 
\lim_{t\nearrow 1}d\pi^0(u)\pi^0(a_t)v_{\hbox{\fiverm K}}\cr 
&=\lim_{t\nearrow 1}\pi^0(a_t) (d\pi^0(\Ad (a_t^{-1})u) 
v_{\hbox{\fiverm K}}\cr 
&=\lim_{t\nearrow 1}\pi^0(a_t) (d\pi^0(\Ad (x_{\hbox {\fiverm H}} 
)^{-1}u)v_{\hbox{\fiverm K}}\cr 
&= c(\Ad(x_{\hbox {\fiverm H}}^{-1})u)v_{\hbox{\fiverm H}}\ .\cr}$$ 
This completes the proof of  the lemma.\qed 
 
\subheadline{2.4. Principal series representations} 
 
\noindent 
In this section we consider the case where $\pi=\pi_\lambda$ is a {\it spherical principal series representation}. 
In particular we will be discuss  the dependence of $v_{\hbox{\fiverm H}}$ 
on the spectral parameter $\lambda$. 
 
\par Let us first recall some well known  facts about the principal series representations. 
For $\alpha\in \Sigma$ let $m_\alpha=\dim \g^\alpha$ and 
$\rho\:={1\over 2}\sum_{\alpha\in \Sigma^+} m_\alpha \alpha$. Write 
$M=Z_K(\a)$ and denote by $\kappa\: G\to K$ and $a\: G\to A$ the projections onto $K$, resp. $A$, 
associated to the Iwasawa decomposition $G=NAK$. Note that $a$ and $\kappa$ have 
unique holomorphic extension 
to $\tilde\Xi$ also denoted by $a$ and $\kappa$ (cf.\ (1.2.1) and [KS01a]). As we are assuming that 
$G\subseteq G_\C$ with $G_\C$ simply connected and $H=G^\tau$, we have $M=Z_H(A)$. In particular 
$M\subseteq H\cap K$. 
 
\par Define a minimal parabolic subgroup of $G$ by $P_{\rm min}=MAN$. 
For $\lambda\in \a_\C^*$ let 
$${\cal D}_\lambda=\{f\in C^\infty(G)\: (\forall g\in G)(\forall man\in P) 
\ f(mang)=a^{\rho-\lambda} f(g)\}\ .$$ 
The group $G$ acts on ${\cal D}_\lambda$ by right translation. Denote the 
corresponding representation by $\pi_\lambda^\infty$, i.e., $\pi_\lambda^\infty(g)f(x)=f(xg)$. 
Denote by $\tilde{\cal H}_\lambda$ the completion of $\cal D_\lambda$ in the norm 
corresponding to the inner product 
$$\la f,g\ra =\int_K f(k)\overline{g(k)}\, dk\, .$$ 
Then $\pi_\lambda^\infty$ extends to a representation $\pi_\lambda$  of $G$ in ${\cal H}_\lambda$.  
We refer to $(\tilde {\cal H}_\lambda ,\pi_\lambda)$ as the {\it spherical principal series representation of $G$ 
with parameter $\lambda$}. The principal series representations $(\pi_\lambda, \tilde {\cal H}_\lambda)$ 
are admissible and of finite length; they are irreducible for generic 
parameters $\lambda$ and unitary if $\lambda\in i\a^*$. It is also well known that 
$\tilde{\cal H}_\lambda^\infty =\cal{D}_\lambda$ excusing our above notation $\pi_\lambda^\infty$. 
 
\par The restriction map $\tilde{\cal{H}}_\lambda\ni f\mapsto f\res_{K}\in L^2(K)$ is 
injective by the left $NA$-covariance of $f$. Furthermore $f\res_{K}$ is left $M$-invariant. Hence 
the restriction map defines an isometry $\tilde{\cal H}_\lambda\hookrightarrow L^2(M\bs K)$. On 
the other hand if $F\in L^2(M\bs K)$ then we can define $f\in  \tilde{\cal H}_\lambda$ by 
$f(nak)=a^{\rho -\lambda}F(Mk)$. Hence 
$\tilde{\cal H}_\lambda \simeq L^2(M\bs K)$. In this realization we 
have 
$$[\pi_\lambda(g)f](Mk)=a(kg)^{\rho -\lambda}f(M\kappa (kg))\, .$$ 
Hence the Hilbert space is the same for all $\lambda$ but the formula for the representation 
depends on $\lambda$. 
Notice that for $k\in K$ this simplifies to 
$(\pi_\lambda(k)f)(Mx)= f(Mxk)$.  We write ${\cal H}_\lambda =L^2(M\bs K)$ to indicate 
the role of $L^2(M\bs K)$ as the representation space for $\pi_\lambda$ and call it 
the {\it compact realization of} $\pi_\lambda$. 
As a consequence of this discussion, we see that $v_{{\hbox{\fiverm K}, \scriptscriptstyle{\lambda}}}=\1_{M\bs K}$ is a 
normalized $K$-fixed vector in 
${\cal H}_\lambda$ and in fact ${\cal H}_\lambda^K=\C v_{{\hbox{\fiverm K}, \scriptscriptstyle{\lambda}}}$. 
We write $v_{{\hbox{\fiverm H}, \scriptscriptstyle{\lambda}}}$ 
and $v_{{\hbox{\fiverm H}, \scriptscriptstyle{\lambda}}}^\omega$ 
instead of $v_{\hbox{\fiverm H}}$ and $v_{\hbox{\fiverm H}}^\omega$ to 
indicate the dependence of $\lambda$. 
 
\par We have $\pi_\lambda^0=\pi_{-\oline \lambda}$. Thus in the compact realization we have 
 
$$[\pi_\lambda^0(a_t)v_{{\hbox{\fiverm K}, \scriptscriptstyle{\lambda}}}](Mk)=a(ka_t)^{\rho+\oline \lambda}\leqno(2.4.1)$$ 
for all $0\leq t<1$ and so 
$$\la v_{{\hbox{\fiverm H}, \scriptscriptstyle{\lambda}}}^\omega, v\ra =\lim_{t\nearrow 1} \int_{M\bs K} a(ka_t)^{\rho+\oline \lambda} 
\, \oline{ v(Mk)} \ dMk \qquad (v\in {\cal H}_\lambda^\omega)\ .\leqno(2.4.2)$$ 
\par In the sequel it will be importnat  
that ${\cal H}_\lambda^\omega=C^\omega(M\bs K)$ is independent of $\lambda$.  
The following theorem specifies the dependence of $v_{{\hbox{\fiverm H}, \scriptscriptstyle{\lambda}}}$  on $\lambda$: 
 
\Theorem 2.4.1. The mapping  
$$\a_\C^*\to \coprod_{\lambda\in \a_\C^*} ({\cal H}_\lambda^{-\infty})^H, 
\ \ \lambda\mapsto v_{{\hbox{\fiverm H}, \scriptscriptstyle{\lambda}}}$$ 
is weakly anti-holomorphic in the sense that for all $v\in C^\omega(M\bs K)$ the mapping  
$$\a_\C^*\to\C, \ \ \lambda\mapsto \la v_{{\hbox{\fiverm H}, \scriptscriptstyle{\lambda}}}, v\ra $$ 
is anti-holomorphic.

\Proof. Let $v\in C^\omega(M\bs K)$. It is convenient to consider  
$v$ as an $M$-invariant function on $K$. As $v$ is analytic, there exists 
an open $K\times K$-invariant neighborhood ${\cal U}$ of $K$ in $K_\C$ such that  
$v$ extends to a holomorphic $M$-invariant function $\tilde v$ on  ${\cal U}$.  
 
\par Let $\eps>0$. Then it follows from the compactness of $Ka_\eps$ and (1.2) that  
we can choose $\eps>0$ small enough such that $\kappa(Ka_\eps)\subeq {\cal U}$.  
 
\par Now consider $v$ as an element of ${\cal H}_\lambda^\omega$. We claim  
that $\pi_\lambda(a_\eps)v$ exists. In fact, using our introductory remarks,  we have   
 
$$[\pi_\lambda(a_\eps)v](Mk)=a(ka_\eps)^{\rho-\lambda} \tilde v(\kappa(ka_\eps))\ .$$ 
 
\par With (2.4.1) we now compute  
 
$$\eqalign{\la v_{{\hbox{\fiverm H}, \scriptscriptstyle{\lambda}}}, v\ra&=\la \pi_\lambda(a_{1-\eps})^0v_{{\hbox{\fiverm K}, \scriptscriptstyle{\lambda}}}, \pi_\lambda(a_\eps)v\ra \cr  
&=\int_{M\bs K} a(ka_{1-\eps})^{\rho+\oline \lambda}\cdot  
 \oline {a(ka_\eps)^{\rho-\lambda}}\cdot \oline {\tilde v(\kappa(ka_\eps))}\ dMk\ .\cr}\leqno(2.4.3)$$ 
By our remarks at the beginning of the proof, we have  
$$\sup_{k\in K} |\tilde v(\kappa(ka_\eps))|<\infty\ .\leqno(2.4.4)  $$ 
Notice that (1.2) implies that both $a(Ka_{1-\eps})$  and $a(Ka_\eps)$ are compact  
subsets of $T(\Omega)$. Thus, if  $C\subeq \a_\C^*$ is  a compact subset, then 
 
$$\sup_{\lambda\in C}\sup_{k\in K} |a(ka_{1-\eps})^{\rho+\oline \lambda}|<\infty,  
\qquad\hbox{and}\qquad   \sup_{\lambda\in C}\sup_{k\in K} |a(ka_\eps)^{\rho -\lambda}|<\infty 
\ .\leqno(2.4.5)$$ 
Therefore, if we use the estimates (2.4.4) and (2.4.5), it follows from (2.4.3) that  
$\lambda\mapsto \la v_{{\hbox{\fiverm H}, \scriptscriptstyle{\lambda}}}, v\ra $ 
is anti-holomorphic.\qed

\subheadline{2.5. Integral representation and asymptotic behaviour of $\theta_\pi$} 
 
\noindent 
Previously we have defined a $H$-invariant holomorphic function  
$\theta_\pi$ for unitary representations $\pi$. For non-unitary $\pi$ we define 
$\theta_\pi$ by  
 
$$\theta_\pi(xK_\C)=\la \pi(x)v_{\hbox{\fiverm K}}, v_{\hbox{\fiverm H}}\ra \qquad (xK_\C\in \Xi_H)\ .$$ 
Clearly, $\theta_\pi$ is a holomorphic function on $\Xi_H$. Moreover, if $\pi$ is unitary, then  
$\theta_\pi$ is  $H$-invariant.

\par In this subsection we will give an integral representation of  
the functions $\theta_\pi$ for principal series represntations $\pi$.  
This  will also allow us to read off the asymptotic behaviour of $\theta_\pi$.  
 
\ssk Recall the definition of the spherical function $\phi_\lambda$ of parameter  
$\lambda\in\a_\C^*$ by  
 
$$\phi_\lambda(g)=\la \pi_\lambda(g)v_{{\hbox{\fiverm K}, \scriptscriptstyle{\lambda}}},  
v_{{\hbox{\fiverm K}, \scriptscriptstyle{\lambda}}}\ra=\int_K a(kg)^{\rho-\lambda}\ dk\qquad (g\in G)\ .$$ 
 
It follows from Theorem 2.1.2 that $\phi_\lambda$ admits a holomorphic extension 
to $\tilde \Xi_H$ (or $\Xi_H$ if we wish to consider $\phi_\lambda$ as a function on $G/K$). 
Also $\phi_\lambda\res_A$ extends holomorphically to the tube $T(2\Omega)=A\exp(2i\Omega)$. 
Notice that 
$z_{\hbox{\fiverm H}}\in T(2\Omega)$. All 
mentioned holomorphic extensions of $\phi_\lambda$ are also denoted by $\phi_\lambda$.

\par In the sequel we abbreviate and write $\theta_\lambda$ instead 
of $\theta_{\pi_\lambda}$. The next result is  immediate 
from the definitions,  Theorem 2.1.2 and the formula (2.4.1). 
 
\Theorem 2.5.1.  Let $(\pi_\lambda, {\cal H}_\lambda)$ be a principal series representation  
with parameter $\lambda\in \a_\C^*$. Then for all $xK_\C\in \Xi_H$ we 
have 
$$\eqalign{ 
\theta_\lambda(xK_\C)&=\lim_{t\nearrow 1} \int_K a(kx)^{\rho-\lambda} \oline{a(ka_t)^{\rho+\oline \lambda}} \ dk\cr 
&=\lim_{t \nearrow 1}\varphi_\lambda (a_tx)\, .}$$ 
Furthermore, 
$$\theta_\lambda(aK_\C)=\phi_\lambda(x_{\hfH}a)$$ 
for all $a\in T(\Omega)=A\exp (i\Omega)$. Here 
$\phi_\lambda$ denotes the holomorphically extended spherical function to $T(2\Omega)$. 
 
\Proof. Let $x\in\Xi_H$. Then we have 
$$\eqalign{ 
\theta_\lambda (x)&=\la \pi_\lambda (x)v_{\hfK},v_{\hfH}\ra\cr 
&=\lim_{t\nearrow 1}\la \pi_\lambda (x)v_{\hfK},\pi^0_{\lambda} (a_t)v_{\hfK}\ra\cr 
&=\lim_{t\nearrow 1}\int_K a (kx)^{\rho -\lambda}\, \overline{a(ka_t)^{\rho +\overline{\lambda}}}\, dk\, .\cr }$$ 
But we can also write the third line as 
$$\eqalign{ 
\theta_\lambda (x)&=\lim_{t\nearrow 1}\la \pi_\lambda (x)v_{\hfK},\pi^0_{\lambda}(a_t)v_{\hfK}\ra \cr 
&=\lim_{t\nearrow 1}\la \pi_\lambda(a_tx)v_{\hfK},v_{\hfK}\ra\cr 
&=\lim_{t\nearrow 1}\varphi_\lambda (a_t x) 
\, .\cr}$$ 
The last statement follows now from Theorem 2.1.2, part (ii). 
\qed 
To discuss the asymptotic expansions of $\theta_\lambda$ along a 
positive  Weyl chamber we first have to recall 
some facts on the Harish-Chandra expansion of the spherical functions. 
For that let $\a_+=\{ X\in \a\: (\forall \alpha\in \Sigma^+) \ \alpha(X)>0\}$ and set 
$A^+=\exp(\a_+)$. Further we define $\Lambda=\N_0[\Sigma^+]$. 
If $\mu\in \Lambda$, then we define a meromorphic function 
$\Gamma_\mu(\lambda)$ in the parameter $\lambda\in \a_\C^*$ by 
$\Gamma_0(\lambda)=1$ and then recursively by 
 
$$\Gamma_\mu(\lambda)={2\over \la \mu, \mu-\lambda\ra} \sum_{\alpha\in \Sigma^+} 
m_\alpha\sum_{k\in \N} \Gamma_{\mu-2k\alpha} \la \mu+\rho-2k\alpha-\lambda,\alpha\ra\ .$$ 
We call $\lambda\in \a_\C^*$ generic if $\Gamma_\mu(\cdot)$ is holomorphic 
at $\lambda$ for all $\mu\in \Lambda$. 
For generic $\lambda\in \a_\C^*$ we define the Harish-Chandra $\Phi$-function 
on $A^+$ by 
$$\Phi_\lambda(a)=a^{\lambda-\rho} \sum_{\mu\in \Lambda} \Gamma_\mu(\lambda) a^{-\mu} \qquad (a\in A^+)\ .$$ 
This series is locally absolutely convergent. In particular, we see that 
$\Phi_\lambda$ extends to a holomorphic function on 
$A^+\exp(2i\Omega)\subeq A_\C$ which we also denote by $\Phi_\lambda$. 
Finally, with $\c(\lambda)$ the familiar Harish-Chandra $\c$-function on $G/K$, we have 
for all generic parameters $\lambda\in \a_\C^*$ that 
$$\phi_\lambda(a)=\sum_{w\in {\cal W}} \c(w\lambda)\Phi_ {w\lambda}(a)\qquad (a\in A^+)\ .$$ 
Combining these facts with Theorem 2.5.1 we now obtain that: 
 
\Theorem 2.5.2. Let $\lambda\in \a_\C^*$ be a generic parameter. Then the following 
assertions hold. 
 
\item{(i)} For $a\in A^+\exp(i\Omega)$ we have 
$$\theta_\lambda(aK_\C)=\phi_\lambda(z_{\hbox{\fiverm H}}a)=\sum_{w\in W} \c(w\lambda) (z_{\hbox{\fiverm H}}a)^{w\lambda 
-\rho}\Phi_{w\lambda}(z_{\hbox{\fiverm H}}a)\ .$$ 
\item{(ii)}  Suppose that $\la \Re\lambda,\alpha\ra>0$ for all $\alpha\in \Sigma^+$. 
Fix $Y\in\a_+$.  Then  
$$\lim_{t\to \infty} e^{t(\rho -\lambda)(Y)} \theta_\lambda(\exp(tY)K_\C) = 
\c(\lambda)\cdot z_{\hbox{\fiverm H}}^{\lambda -\rho} \ .$$\qed

\subheadline{2.6. $H$-orbit coefficients of the holomorphic distribution vector} 
 
\noindent 
As we have remarked already earlier the space $({\cal H}_\lambda^{-\infty })^H$ 
has dimension $\left |{\cal W}/ {\cal W}_0\right|$ for generic $\lambda$. 
One can parametrize  $({\cal H}_\lambda^{-\infty })^H$ through the open 
$H$-orbits in the flag manifold $P_{\rm min}\bs G$. These orbits haven been parametrized 
by Rossmann and Matsuki (cf. [M79]); they are given by 
$$P_{\rm min} wH \qquad (w\in  {\cal W}/ {\cal W}_0)\ .$$ 
 
\par For $\lambda\in \a_\C^*$ and $w\in  {\cal W}/ {\cal W}_0$ define 
a right $H$-invariant function on $G$ by 
 
$$\eta_{\lambda, w}(x)=\cases{ a^{\rho+\oline \lambda} & for $x=manwh\in MAN wH$ \cr 
0 & otherwise\ .\cr}$$ 
 
For $\lambda\in \a^*$ we will use the notation $\lambda <<0$ if 
$\la \lambda,\alpha\ra<<0$ for all $\alpha\in \Sigma^+$. 
Then it is known that for $\lambda<<0$ the functions $\eta_{\lambda, w}$ are 
continuous and define $H$-fixed distribution vectors of $\pi_\lambda$ 
[\'O87]. 
Moreover,  the distributions $\eta_{\lambda, w}$ admit continuation 
in $\lambda$ to a weakly anti-meromorphic function on $\a_\C^*$. For generic $\lambda\in \a_\C^*$ 
we have 
$$({\cal H}_\lambda^{-\infty})^H=\bigoplus_{w\in {\cal W}/{\cal W}_0} \C \eta_{\lambda,w}$$ 
and the mapping 
$$j_\lambda\: \C^{\left |{\cal W}/{\cal W}_0\right|}\to ({\cal H}_\lambda^{-\infty})^H, 
\ \ (c_w)_w\mapsto \sum_{w\in {\cal W}/{\cal W}_0} c_w \eta_{\lambda,w}$$ 
is a bijection for generic $\lambda$, weakly anti-meromorphic in $\lambda$ [vdB88]. 
For $\lambda<<0$ the inverse of $j_\lambda$ 
is given by the evaluation mapping 
$${\rm ev}\: ({\cal H}_\lambda^{-\infty})^H\to \C^{\left |{\cal W}/{\cal W}_0\right|}, 
\ \ \eta\mapsto (\eta(w))_w\ .$$ 
 
\par On the other hand we know that the weakly holomorphic distribution 
vector $v_{{\hbox{\fiverm H}, \scriptscriptstyle{\lambda}}}$ 
depends weakly anti-holomorphically on  $\lambda\in \a_\C^*$ (cf.\ Theorem 2.4.1). 
The next theorem gives us the coefficients of $v_{{\hbox{\fiverm H}, \scriptscriptstyle{\lambda}}}$ 
in terms of the basis $(\eta_{\lambda,w})$ of $({\cal H}_\lambda^{-\infty})^H$. 
We note that for $\lambda<<0$ the distribution $v_{{\hbox{\fiverm H}, \scriptscriptstyle{\lambda}}}$ is given through 
the bounded measurable function 
 
$$v_{{\hbox{\fiverm H}, \scriptscriptstyle{\lambda}}}(k)=\lim_{t\nearrow 1} a(ka_t)^{\rho+\oline \lambda}\ .\leqno(2.6.1)$$

\Theorem 2.6.1. For generic parameters $\lambda\in \a_\C^*$ we have 
 
$$v_{{\hbox{\fiverm H}, \scriptscriptstyle{\lambda}}}=\sum_{w\in {\cal W}/ {\cal W}_0} 
 z_{\hbox{\fiverm H}}^{w^{-1}(\rho +\oline\lambda)}\cdot \eta_{\lambda,w}\ .$$ 
 
\Proof. For $\lambda\in \a_\C^*$ generic write 
$$v_{{\hbox{\fiverm H}, \scriptscriptstyle{\lambda}}}=\sum_{w\in {\cal W}/ {\cal W}_0} c_{\lambda, w} \cdot  \eta_{\lambda,w}$$ 
for the basis expansion. As  the coefficients $  c_{\lambda, w}$ depend weakly 
anti-meromorphically on $\lambda$, it is sufficient to show that 
$c_{\lambda, w}=z_{\hbox{\fiverm H}}^{w^{-1}(\rho 
+\oline \lambda)}$ 
for $\lambda<<0$. Then (2.6.1) implies that 
 
$$c_{\lambda,w}=v_{\hfH}(w)=\lim_{t\nearrow 1} v_{{\hbox{\fiverm H}, \scriptscriptstyle{\lambda}}}(wa_t)=a(wz_{\hbox{\fiverm H}})^{\rho+\oline \lambda}= 
z_{\hbox{\fiverm H}}^{w^{-1}(\rho+\oline \lambda)}\ , $$ 
as was to be shown. \qed 
 
\Remark 2.6.2. Let us go back to Example 2.1.5 for $G=\Sl(2,\R)$. Here 
${\cal W}_0=\{\1\}$ and so ${\cal W}/{\cal W}_0=\{ \1, w\}$ where $w$ is the non-trivial element in 
the Weyl group which acts by multiplication by $-1$. 
The distributions $v_{\hbox{\fiverm H,1}}$ and  $v_{\hbox{\fiverm H,2}}$ from Example 2.1.5 are given in 
the above notation by

$$v_{\hbox{\fiverm H,1}}=\eta_{\lambda,\1}\qquad\hbox{and}\qquad v_{\hbox{\fiverm H,2}}=\eta_{\lambda,w}\ .$$ 
In Example 2.1.5 we did show that 
$$v_{\hbox{\fiverm H}}=v_{{\hbox{\fiverm H}, \scriptscriptstyle{\lambda}}}=e^{i{\pi\over 4} 
(1+\overline{\lambda})}v_{\hbox{\fiverm H,1}}+ e^{-i{\pi\over 4}(1+\overline{\lambda})} v_{\hbox{\fiverm H,2}}\. $$ 
As $z_{\hbox{\fiverm H}}^{\rho+\overline{\lambda}}=e^{i{\pi\over 4}(1+\overline{\lambda})}$, we hence see that 
the above formula is a special case of Theorem 2.6.1.\qed

\subheadline{2.7. Transformation under the intertwining matrix} 
 
\noindent 
For $w\in {\cal W}$ and $\lambda\in \a_\C^*$ generic  we have an intertwining operator 
$$A(\lambda, w\lambda)\: (\pi_\lambda, {\cal H}_\lambda^\infty)\to (\pi_{w\lambda}, 
{\cal H}_{w\lambda}^\infty)\ .$$ 
Notice that the Hilbert space ${\cal H}_\lambda$ is independent of $\lambda$ because 
we use the compact realzation. We can therefore speak about meromorphic maps 
from $\a_\C^*$ into the space of bounded operators from 
${\cal H}_\lambda$ into ${\cal H}_\mu$.  In this sense 
it is well known that the map $\a_\C^*\ni \lambda\mapsto  A(\lambda, w\lambda)$ is meromorphic. 
Dualizing, we obtain an anti-meromorphic family of intertwining operators 
$$A(\lambda, w\lambda)^*\: {\cal H}_{w\lambda}^{-\infty}\to {\cal H}_\lambda^{-\infty}\ . $$ 
Restricting $A(\lambda ,w\lambda)^*$ to the space of $H$-invariant distribution 
vectors we obtain a linear bijection, say 
 
$$A_H(\lambda, w\lambda)^*\: ({\cal H}_{w\lambda}^{-\infty})^H\to ({\cal H}_\lambda^{-\infty})^H\ .$$ 
Often one refers to  $A_H^*(\lambda, w\lambda)$ as the {\it intertwining  matrix}. 
In terms of the basic distribution vectors $(\eta_{w\lambda, w'})_{w'}$ respectively 
$(\eta_{\lambda, w'})_{w'}$ 
the operator  $A_H(\lambda, w\lambda)^*$ has an unknown, seemingly complicated expression. 
In this section we will show that the the intertwining  matrix 
maps the holomorphic distribution vector 
$v_{{\hbox{\fiverm H}, \scriptscriptstyle{w\lambda}}}$ 
to a multiple of $v_{{\hbox{\fiverm H}, \scriptscriptstyle{\lambda}}}$. 
In order to describe this multiple more precisely we need more notations. 
 
\par For $w\in {\cal W}$ define a subgroup of $\oline N=\theta(N)$ by 
$$\oline N_w=\oline N\cap wNw^{-1}\ .$$ 
For $\Re \lambda<<0$ we define functions 
$$\c_w(\lambda)=\int_{\oline N_w} a(\oline n)^{\rho-\lambda} \ d\oline n\ .$$ 
If $w=w_0$ is the longest element in ${\cal W}$, then we write 
$\c(\lambda)$ instead of $\c_{w_0}(\lambda)$ and remark that 
$\c(\lambda)$ is the familiar Harish-Chandra $c$-function on $G/K$. 
The functions $\c_w(\lambda)$ admit meromorphic continuation to $\a_\C^*$ 
and can be explicitely computed (Gindikin-Karpelevic formula). 
 
\par With this notation the intertwinig operators $A(\lambda,w\lambda)$  for $\lambda<<0$ 
are defined by 
$$[A(\lambda, w\lambda)f](x) =\int_{\oline N_w} f(\oline n w x) \ d\oline n 
\qquad (f\in {\cal D}_\lambda, x\in G)\ .$$

\Theorem 2.7.1. Assume that $\lambda\in \a_\C^*$ is generic. Then 
$$A_H(\lambda, w\lambda)^*v_{{\hbox{\fiverm H}, \scriptscriptstyle{w\lambda}}} 
= \c_w(\oline\lambda)\cdot v_{{\hbox{\fiverm H}, \scriptscriptstyle{\lambda}}}\ .$$

\Proof. It is well known - and follows immediately from 
the definition -  that  $A(\lambda, w\lambda)v_{{\hbox{\fiverm K}, \scriptscriptstyle{\lambda}}}=\c_w(\lambda) 
v_{{\hfK,w\lambda}}$. But then, as $A(\lambda, w\lambda )^*v_{\hfK,w\lambda}$ is $K$-invariant, we also get 
$$\eqalign{ 
\la A(\lambda ,w\lambda )^*v_{\hfK,w\lambda},v_{\hfK,\lambda}\ra &= \la v_{\hfK,w\lambda},A(\lambda,w\lambda)v_{\hfK,\lambda}\ra\cr 
&=\la v_{\hfK,w\lambda },\c_w(\lambda )v_{\hfK,w\lambda}\ra\cr 
&=\overline{\c_w (\lambda )}\, .\cr} 
$$ 
Noticing that $\overline{\c_w(\lambda)}=\c_w(\overline{\lambda})$ it follows that 
$$A(\lambda ,w\lambda )^*v_{\hfK,w\lambda}=\c_w(\overline{\lambda})v_{\hfK,\lambda}\, .$$ 
Finally, using that $A(\lambda ,w\lambda )^*$ is an intertwining operator, we get for 
a $K$-finite  $u$ : 
$$\eqalign{ 
\la A_H^*(\lambda ,w\lambda)v_{\hfH,w\lambda},u\ra &= 
\la v_{\hfH,w\lambda},A(\lambda ,w\lambda )u\ra\cr 
&=\lim_{t\nearrow 1}\la \pi_{w\lambda }^0(a_t)v_{\hfK,w\lambda},A(\lambda,w\lambda )u\ra\cr 
&=\lim_{t\nearrow 1}\la v_{\hfK,w\lambda},\pi_{w\lambda }(a_t)^{-1}A(\lambda,w\lambda )u\ra\cr 
&=\lim_{t\nearrow 1}\la v_{\hfK,w\lambda},A(\lambda,w\lambda )\pi_{w\lambda }(a_t)^{-1}u\ra\cr 
&=\lim_{t\nearrow 1}\la A(\lambda,w\lambda )^*v_{\hfK,w\lambda},\pi_{w\lambda }(a_t)^{-1}u\ra\cr 
&= \c_w(\overline{\lambda})\lim_{t\nearrow 1}\la \pi_{\lambda }(a_t)v_{\hfK,\lambda },u\ra\cr 
&=\c_w (\overline{\lambda}) \la v_{\hfH,\lambda },u\ra \, .\cr} 
$$ 
Hence $A_H(\lambda,w\lambda)^*v_{{\hbox{\fiverm H}, \scriptscriptstyle{w\lambda}}}= \c_w(\oline \lambda)\cdot 
v_{{\hbox{\fiverm H}, \scriptscriptstyle{\lambda}}}$, as was to be shown. \qed

\subheadline{2.8. Relation to the horospherical picture} 
 
\noindent 
In this subsection we explain the construction of the 
holomorphic distribution vector $v_{\hbox{\fiverm H}}$ from 
the horospherical point of view. 
 
\par Let $\kappa\in \a_\C^*$ be a complex parameter. 
We call the holomorphic function 
$$ \tilde \psi _\kappa(nak)= a^\kappa=\exp(< \kappa , \log a>), 
\qquad n\in N_{\Bbb C}, a\in T(\Omega), k\in K_\Bbb C$$ on $\tilde 
\Xi $ (see\ (1.2.1)) the {\it holomorphic horospherical function 
with parameter $\kappa$}. 
 
\par This function can be pushed down to $\Xi$ as a holomorphic locally $N_\C$-invariant function 
$\psi_\kappa$. The function $\psi_\kappa$ is a holomorphic extension of the usual 
horospherical ($N$-invariant, $A$-homogeneous) function on $G/K$, 
corresponding to spherical principal representations related to $\kappa$.  
For certain  values of the parameter $\kappa$ the function $\psi_\kappa$ has boundary distribution values 
$\psi_{\kappa, \hbox{\fiverm H}}$ on $G/H$. To 
understand the structure of those distributions, let us remark that a Zariski open part 
of G/H is the disjoint union of domains $Y_j= NAy_j$ where $y_1, 
\cdots, y_k$ correspond to the vertices of $\Omega_H$ 
(${\cal W}$-equivalent). On each $Y_j$ we have an $N$-invariant, $A$-homogeneous  
distribution with parameter $\kappa$ (such 
distributions on $Y_j$ are unique up to a multiplicative 
constant). 
 
The function $\tilde \psi_\kappa$  can be also pushed down on a domain 
$D_H$ in $N_\C M_\C \bs G_\C$ as a holomorphic function 
$\psi_\kappa$. This function with parameters holomorphically 
extends $K$-invariant vectors in the principal spherical 
representations on $MN\bs G $. The boundary values of this holomorphic 
function give an $H$-invariant distribution, the domains $Y_j$ 
correspond to the $H$-orbits on $MN\bs G$  and we have the corresponding 
decomposition of $\psi_\kappa$.

\sectionheadline{3. An application: Hardy spaces for NCC symmetric spaces} 
 
\noindent 
In this section we apply  our theory 
developed in Section 2 to associate 
to every NCC symmetric space $G/H$ a Hardy space ${\cal H}^2(\Xi_H)$. 
The Hardy space is a $G$-invariant Hilbert space 
of holomorphic functions on $\Xi_H$ featuring a boundary value mapping which gives 
an isometric embedding of  ${\cal H}^2(\Xi_H)$ into the most-continuous 
spectrum $L^2(G/H)_{\rm mc}$ of $L^2(G/H)$. Hence we give a 
realization of a part of $L^2(G/H)_{\rm mc}$ in a space 
of holomorphic function on $\Xi_H$, generalizing and extending our 
previous results from [GK\'O01] to all NCC spaces. 
 
\par This section is organized as follows: After a brief digression on 
Hilbert spaces of holomorphic functions on $\Xi_H$, we  give 
an adhoc definition of the Hardy space through the spectral measure. Then, 
after recalling the theory of the most-continuous spectrum, we will show that 
there is a boundary value mapping embedding ${\cal H}^2(\Xi_H)$ isometrically 
into $L^2(G/H)_{\rm mc}$.

\subheadline{3.1. $G$-invariant Hilbert spaces of holomorphic functions on $\Xi_H$} 
 
\noindent 
In this section we briefly recall the abstract theory of $G$-invariant 
Hilbert spaces of holomorphic functions, specialized to the complex manifold $\Xi_H$ 
(see  also [FT99] and [K99] for the general theory). 
 
\par In the sequel we will consider ${\cal O}(\Xi_H)$ as a Fr\'echet space with the 
topology of compact convergence. We let $G$ act on ${\cal O}(\Xi_H)$ by the 
left regular representation $L$: 
$$(L(g)f)(z)=f(g^{-1}z)\qquad (g\in G, f\in {\cal O}(\Xi_H), z\in \Xi_H)\ .\leqno(3.1.1)$$ 
By a {\it $G$-invariant Hilbert space 
of holomorphic functions on $\Xi_H$} we understand a Hilbert space ${\cal H}\subeq {\cal O}(\Xi_H)$ 
such that:

\ssk 
\item{(IH1)} The inclusion ${\cal H}\into {\cal O}(\Xi_H)$ is continuous. 
\item{(IH2)} The Hilbert space ${\cal H}$ is invariant under $L$ and the the 
corresponding representation of $G$ is unitary. 
\ssk 
It follows from (IH1) that for every $z\in \Xi_H$ the  point 
evaluation ${\cal H}\to \C, \ f\mapsto f(z)$, 
is continuous. Thus, there exists a ${\cal K}_z\in {\cal H}$ such that 
$\la f, {\cal K}_z\ra =f(z)$ holds for every $f\in {\cal H}$. In 
this way we obtain a function 
$${\cal K}\: \Xi_H \times \Xi_H\to \C, \ \ (z,w)\mapsto 
{\cal K}(z,w)=\la {\cal K}_w, {\cal K}_z\ra\ .$$ 
The function ${\cal K}$  is holomorphic in the first variable and anti-holomorphic in the second variable. 
It follows from (IH2) that ${\cal K}$ is $G$-invariant, i.e., ${\cal K}(gz,gw)= 
{\cal K}(z,w)$ holds for all $g\in G$ and all $z,w\in \Xi_H$. 
We call ${\cal K}$ the {\it Cauchy-Szeg\"o kernel} of ${\cal H}$ and note that 
${\cal H}$ is determined by ${\cal K}$. 
 
\msk  To describe the spectral resolution of ${\cal K}$ denote by $\hat G_s$  the $K$-spherical unitary dual of $G$. 
We view $\hat G_s$ as a subset of $\a_\C^*/ {\cal W}$ using the 
parametrization of the spherical principal series. Notice that the topology on 
$\hat G_s$ coincides with the topology induced from  
$\a_\C^*/{\cal W}$. 
Slightly abusing our notation from Subsection 2.4, we denote by 
$(\pi_\lambda, {\cal H}_\lambda)$ a representative of $\lambda\in \hat G_s$.  
For $\lambda\in \hat G_s$ define the $G$-invariant kernel ${\cal K}_\lambda$ by 
$${\cal K}_\lambda(xK_\C, yK_\C)=\la \pi_\lambda(\oline y)v_{{\hbox{\fiverm K}, \scriptscriptstyle{\lambda}}}, 
 \pi_\lambda(\oline x)v_{{\hbox{\fiverm K}, \scriptscriptstyle{\lambda}}}\ra 
\qquad (xK_\C, yK_\C\in \Xi_H)\ .$$ 
Write $\Xi_H^{\rm opp}$ for $\Xi_H$ but endowed with the  
opposite complex structure. We recall that the map  
$$\hat G_s \to {\cal O}(\Xi_H\times \Xi_H^{\rm opp}),  
\ \ \lambda\mapsto {\cal K}_\lambda$$ 
is continuous [KS01b, Sect. 5]  
(it follows from the fact that the spherical  
functions $\phi_\lambda$ and their holomorphic continuations   
are continuous in $\lambda$).  
Then by [KS01b, Th.\ 5.1] there exists a unique Borel measure $\mu$ on 
$\hat G_s$ such that 
 
$${\cal K}(z,w)=\int_{\hat G_s} {\cal K}_\lambda(z,w)\ d\mu(\lambda) \qquad (z,w\in \Xi_H)\leqno(3.1.2)$$ 
with the  right hand side converging absolutely on compact subsets of $\Xi_H\times \Xi_H$. 
Equivalently phrased, the mapping 
 
$$\Phi\: \int_{\hat G_s}^\oplus {\cal H}_\lambda\ d\mu(\lambda)\to {\cal H} 
, \ \ s=(s_\lambda)_\lambda\mapsto \left(xK_\C\mapsto \int_{\hat G_s} \la \pi_\lambda(x^{-1})s_\lambda, 
v_{{\hbox{\fiverm K}, \scriptscriptstyle{\lambda}}}\ra \ d\mu(\lambda)\right)\leqno (3.1.3)$$ 
is a $G$-equivariant unitary isomorphism. In the sequel we refer to the measure $\mu$ as 
the {\it Plancherel measure} of ${\cal H}$. 
 
\ssk  In [KS01b] 
a criterion was given on a Borel measure $\mu$ on $\hat G_s$ to be a Plancherel measure 
for an invariant Hilbert space ${\cal H}={\cal H}(\mu)$ on $\Xi$. This criterion can be easily adapted to 
invariant Hilbert spaces on $\Xi_H$. Let us provide the necessary modifications. 
\par Define a norm $\|\cdot\|_{\hbox{\fiverm H}}$ on $\a_\C^*$ by 
 
$$\|\lambda\|_{\hbox{\fiverm H}}\:=\sup_{w\in {\cal W}/{\cal W}_0} |\lambda(wX_{\hbox{\fiverm H}})|\qquad (\lambda\in \a_\C^*) \ .$$ 
Then [KS01b, Prop. 5.4] and its proof readily gives the following generalization: 
 
\Proposition 3.1.1. Let $\mu$ be a Borel measure on $\hat G_s$ with the property 
$$(\forall 0\leq c<2)\qquad \int_{\hat G_s} e^{c\|\Im \lambda\|_{\hbox{\fiverm H}}}\ d\mu(\lambda)<\infty\ .\leqno(3.1.4)$$ 
Then $\mu$ is the  Plancherel measure of an  
invariant Hilbert space ${\cal H}(\mu)$ on $\Xi_H$. \qed 
 
\subheadline{3.2. The definition of the Hardy space} 
 
\noindent 
We are now ready to give the definition of the Hardy space on $\Xi_H$. 
Let us denote by $i\a_+^*$ an open Weyl chamber in  $i\a^*$. In the sequel we will 
consider $i\a_+^*$ mainly as a subset of $\hat G_s$. Let 
$\cal O$ be a neighborhood of $i\a^*$ such that $\sum_{w\in {\cal W}/ {\cal W}_0} z_{\hbox{\fiverm H}}^{-2w^{-1}\lambda}$ 
has a holomorphic square root 
${\bf z}_{\hfH}(\lambda)$ on $\cal O$. Define a holomorphic function $\cH$ on 
$\cal O$ by 
$$\cH (\lambda)=\c(\lambda)\cdot {\bf z}_{\hfH}(\lambda )$$
and define a Borel measure $\mu$ on $i\a_+^*$ by 
$$d\mu (\lambda )=\frac{d\lambda}{|\cH (\lambda )|^2} \leqno(3.2.1)$$
where $d\lambda$ denotes the Lebesgue measure. Then we have:

\Lemma 3.2.1. The measure $\mu$ satisfies the  condition (3.1.4); in particular
$\mu$ is the Plancherel measure of an invariant Hilbert space ${\cal H}(\mu)$ on $\Xi_H$.

\Proof. Recall the growth behaviour of the $\c$-function on the imaginary axis: There
exists constants $C, N>0$ such that
$$(\forall \lambda\in i\a^*)\qquad {1\over |\c(\lambda)|^2} \leq C (1+|\lambda|^N)\ . $$
Moreover for $\lambda \in i\a^*$ one has 
$z_{\hbox{\fiverm H}}^{\lambda}=e^{\lambda (iX_{\hbox{\fiverm H}})}>0$. 
Hence
$$\sum_{w\in {\cal W}/ {\cal W}_0} |z_{\hbox{\fiverm H}}^{-w^{-1}\lambda}|^2\geq
e^{2\|{\rm Im}\lambda\|_{\hbox{\fiverm H}}}\. $$

Combining these two facts now yields that $\mu$ satisfies (3.1.4).\qed

Using Proposition 3.1.1 and Lemma 3.2.1 we now can give an adhoc-definition of the  Hardy 
space on $\Xi_H$. 
 
\Definition 3.2.2. {\bf(Hardy space)} Let $G/H$ be a NCC symmetric space and $\Xi_H$ its associated domain 
in $G_\C/ K_\C$. Then we define the {\it Hardy space} ${\cal H}^2(\Xi_H)$ on $\Xi_H$ by 
$${\cal H}^2(\Xi_H)={\cal H}(\mu)$$ 
with $\mu$ as in (3.2.1). \qed

Recall the Cauchy-Szeg\"o kernel ${\cal K}(z,w)$ of the invariant 
Hilbert space ${\cal H}^2(\Xi_H)$ from Subsection 3.1. 
 
\Lemma 3.2.3. Let ${\cal K}$ be the Cauchy-Szeg\"o kernel of ${\cal H}^2(\Xi_H)$. 
Then the limits 
$$\Psi(z)=\lim_{t\nearrow 1} {\cal K}(z,a_tK_\C)\qquad (z\in \Xi_H)$$ 
exist locally uniformly. In particular, $\Psi\:\Xi_H\to \C$ 
is an $H$-invariant holomorphic function. 
 
\Proof. Fix $z\in \Xi_H$ and let $U\subeq \Xi_H$ be a compact 
neighborhood of $z$. Choose $\eps>0$ small enough such that 
$a_\eps U\subeq \Xi_H$. Then, by $G$-invariance, we have 
$${\cal K}(z, a_tK_\C)={\cal K}(a_\eps z, a_{-\eps}a_tK_\C)= 
{\cal K}(a_\eps z,a_{t-\eps}K_\C)$$ 
for all $\eps <t<1$. The claim follows now, because 
$]\eps, 1+\eps [\ni t\mapsto {\cal K}(a_\eps z,a_{t-\eps}K_\C)\in \C$ 
is continuous, and hence 
$\lim_{t\nearrow 1}{\cal K}(z, a_tK_\C)={\cal K}(a_{\eps} z, a_{1-\eps}K_\C)$ 
exists and the convergence is uniform on compact subsets.\qed 
%
 
We refer to 
$\Psi$ as the {\it Cauchy-Szeg\"o function} of ${\cal H}^2(\Xi_H)$. 
As ${\cal K}$ is $G$-invariant, it follows that ${\cal K}$  
can be reconstructed from $\Psi$.  
Moreover, as $H_\C T(\Omega_{\hbox{\fiverm H}})K_\C/ K_\C$ 
meets $\Xi_H$ in an open set, we conclude that $\Psi$ is 
uniquely determined by its restriction to $T(\Omega_{\hbox{\fiverm H}})K_\C/ K_\C \subeq \Xi_H$. 
 
Using Theorem 2.5.2 (i) we finally obtain the spectral 
resolution of $\Psi$.

\Theorem  3.2.4. For $a\in T(\Omega_{\hbox{\fiverm H}})$ we have 
 
$$\Psi(aK_\C)=\int_{i\a_+^*} \phi_\lambda(z_{\hbox{\fiverm H}}a) 
\ {d\lambda\over |\cH(\lambda)|^2}, $$ 
where the integrals on the right hand side converge uniformly 
and absolutely on compact subsets of $T(\Omega_{\hbox{\fiverm H}})$.\qed

\ssk We now discuss the boundary value map $b\: {\cal H}^2(\Xi_H)\to 
L^2(G/H)_{\rm mc}$. As usual, this boundary value map can be 
nicely defined pointwise only on an appropriate dense subspace of ${\cal H}^2(\Xi_H)$. 
Write ${\cal H}^2(\Xi_H)^\omega$ for the analytic  
vectors of the left regular representation  
$(L, {\cal H}^2(\Xi_H))$. Fix $f\in 
{\cal H}^2(\Xi_H)^\omega$. Then for every  
compact subset $C\subeq G$ there exists  
an $0<\eps<1$ such that $L(a_{-\eps}g^{-1})f$ exists 
for all $g\in G$. In particular, if $0<\eps\leq t<1$, then  
$$f(ga_tK_\C)=f(ga_\eps a_{t-\eps} K_\C)= 
[L(a_{-\eps}g^{-1})f)](a_{t-\eps}K_\C)\ .$$ 
and so  
$$\lim_{t\nearrow 1} f(ga_tK_\C)= 
[L(a_{-\eps}g^{-1})f)](a_{1-\eps}K_\C)\ .$$ 
It follows that we have a well defined  
$G$-equivariant boundary value map: 
 
$$b^\omega\: {\cal H}^2(\Xi_H)^\omega\to C(G/H), 
\ \  b^\omega(f)(gH)= 
\lim_{t\nearrow 1} f(ga_tK_\C)\ \ .\leqno(3.2.2)$$ 
 
\msk Recall from (3.1.3) the isomorphism 
$\Phi\: 
\int_{i\a_+^*}^\oplus{\cal H}_\lambda\ d\mu(\lambda) 
\to {\cal H}^2(\Xi_H)={\cal H}^2(\mu)$: 
$$s=(s_\lambda)_\lambda\mapsto \left(xK_\C\mapsto \int_{\hat G_s} \la \pi_\lambda(x^{-1})s_\lambda, 
v_{{\hbox{\fiverm K}, \scriptscriptstyle{\lambda}}}\ra \ d\mu(\lambda)\right)\, .$$ 
It is 
useful to have the corresponding formula for $b^\omega$ 
on the space of 
sections with values in 
$\left(\int_{i\a_+^*}^\oplus{\cal H}_\lambda\ d\mu(\lambda)\right)^\omega$, i.e., on 
the space of 
analytic sections. 
In this regard, it is better to replace 
${\cal H}^2(\Xi_H)^\omega$ by some smaller  
but dense subspace ${\cal H}^2(\Xi_H)_0$.  
In order to define  ${\cal H}^2(\Xi_H)_0$ 
we have to introduce some terminology.  
For a section $s=(s_\lambda)_\lambda\in  
\int_{i\a_+^*}^\oplus{\cal H}_\lambda\  
d\mu(\lambda)$ we define its support  
by $\supp(s)=\oline{\{\lambda\in i\a_+^*\: s_\lambda 
\neq 0\}}$. Futhermore we shall use  
the identifications ${\cal H}_\lambda^\omega= 
C^\omega(M\bs K)$ for $\lambda\in i\a_+^*$.  
Recall that if $f\in{\cal H}^2(\Xi_H)^\omega$ 
and $s=(s_\lambda)_\lambda=\Phi^{-1}(f)$, then  
almost each stalk $s_\lambda$ is an analytic  
vector, i.e. $s_\lambda\in C^\omega(M\bs K)$.  
The subspace ${\cal H}^2(\Xi_H)_0$ is  
then defined by 
$${\cal H}^2(\Xi_H)_0= 
\left\{ f\in {\cal H}^2(\Xi_H)^\omega\quad : \quad \eqalign{& 
f \quad \hbox {is $K$-finite, }\cr 
& s=(s_\lambda)_\lambda=\Phi^{-1}(f) 
\quad \hbox{has compact support, }\cr 
& s\: i\a_+^*\to C^\omega(M\bs K) 
\quad \hbox{is weakly smooth.}}\right\}$$ 
It is an easy verification that  
${\cal H}^2(\Xi_H)_0$ is a dense subspace 
of ${\cal H}^2(\Xi_H)$. Write $b_0\: {\cal H}^2(\Xi_H)_0\to 
C(G/H)$ for the restriction of $b^\omega$ to ${\cal H}^2 
(\Xi_H)_0$.  
 
\par In the sequel we will often identify a function  
$f\in {\cal H}^2(\Xi_H)_0$ with its corresponding  
section $s=(s_\lambda)=\Phi^{-1}(f)$.  
We then claim that  
 
$$b_0\: {\cal H}^2(\Xi_H)_0\to C(G/H), \ \ s=(s_\lambda) 
\mapsto\left(gH\mapsto\int_{i\a_+^*} 
\la \pi_\lambda(g^{-1})s_\lambda, v_{{\hbox{\fiverm H}, 
\scriptscriptstyle{\lambda}}}\ra \ d\mu(\lambda)\right)\leqno(3.2.3)$$ 
Notice that it is a priori not even clear that the  
right hand side of (3.2.3) is 
well defined. To establish (3.2.3) fix  
$f\in {\cal H}^2(\Xi_H)_0$ and $g\in G$.  
As $f$ is an analytic  
vector for the left regular representation  
$(L, {\cal H}^2(\Xi_H))$ it follows that  
there exists an $0<\eps<1$ such that $L(a_\eps g^{-1})f$  
exists. Using standard procedures one deduces   
that  
$\pi_\lambda(a_\eps g^{-1})s_\lambda$ exists  
for almost all $\lambda$. In particular  
$s_\lambda\in {\cal H}_\lambda$  
is analytic for almost all $\lambda$.  
 Furthermore, 
$L(a_\eps g^{-1})f$ corresponds to the section  
$(\pi_\lambda(a_\eps g^{-1})s_\lambda)_\lambda$ and so 
$$\|L(a_\eps g^{-1})f\|^2=\int_{i\a_+^*} 
\|\pi_\lambda(a_\eps g^{-1})s_\lambda\|^2 \ d\mu(\lambda) 
<\infty\ .\leqno(3.2.4)$$ 
 
With the  
convention $\pi_\lambda(a_1) 
v_{{\hbox{\fiverm K }, \scriptscriptstyle{\lambda}}} 
=v_{{\hbox{\fiverm H }, \scriptscriptstyle{\lambda}}}$ 
we then have for all $\eps\leq t\leq 1$ and almost all $\lambda$  
the estimate 
$$\eqalign{|\la \pi_\lambda(g^{-1})s_\lambda,  
\pi_\lambda(a_t)v_{{\hbox{\fiverm K}, \scriptscriptstyle{\lambda}}}\ra| 
&= |\la \pi_\lambda(a_\eps g^{-1})s_\lambda,  
\pi_\lambda(a_{t-\eps}) v_{{\hbox{\fiverm K}, 
 \scriptscriptstyle{\lambda}}}\ra|\cr  
&\leq  \|\pi_\lambda(a_\eps g^{-1})s_\lambda\| 
\cdot \|\pi_\lambda(a_{t-\eps}) v_{{\hbox{\fiverm K}, 
 \scriptscriptstyle{\lambda}}}\|\cr  
&\leq M\cdot   \|\pi_\lambda(a_\eps g^{-1})s_\lambda\|\cr}\leqno 
(3.2.5)$$ 
with $M=\sup_{\lambda\in \supp(s)\atop  
\eps\leq t\leq 1}  
\|\pi_\lambda(a_{t-\eps}) v_{{\hbox{\fiverm K}, 
 \scriptscriptstyle{\lambda}}}\|<\infty$ as $\supp(s)$ is  
compact.  
\par Recall our notion of holomorphic extension from  
Definition 2.2.1.  As almost  
each stalk $s_\lambda$ is  
an analytic vector in ${\cal H}_\lambda$, it follows  
from  estimates (3.2.4-5) and the compactness  
of $\supp (s)$ that  
$$\eqalign{\int_{i\a_+^*} 
 \la \pi_\lambda(g^{-1})s_\lambda,  
v_{{\hbox{\fiverm H}, \scriptscriptstyle{\lambda}}}\ra \ d\mu(\lambda) 
&= 
\int_{\supp (s)} 
\lim_{t\nearrow 1} \la \pi_\lambda(g^{-1})s_\lambda,  
\pi(a_t)v_{{\hbox{\fiverm K}, \scriptscriptstyle{\lambda}}}\ra \ d\mu(\lambda)\cr 
&= 
\lim_{t\nearrow 1}  
\int_{\supp(s)}\la \pi_\lambda(g^{-1})s_\lambda,  
\pi_\lambda(a_t)v_{{\hbox{\fiverm K}, \scriptscriptstyle{\lambda}}}\ra \ d\mu(\lambda)\cr 
&= b^\omega(f)(gH)\ .\cr}$$ 
As $f$ and $g$ were arbitray, this completes the proof  
of (3.2.3).

\subheadline{3.3. The Plancherel Theorem for $L^2(G/H)_{\rm mc}$} 
 
\noindent 
Before we can  show that $b_0$ has image in $L^2(G/H)_{\rm mc}$ and 
extends to an isometric embedding, we need to recall some facts about the 
most continuous spectrum  $L^2(G/H)_{\rm mc}$ of $L^2(G/H)$  (cf.\ [vdBS97a] and [D98]). 
All the results collected below are proved in [vdBS97a] or might be considered as special 
cases of [D98]. The crucial way where our assumption that $G/H$ 
is NCC, $H=G^\tau$, and $G\subseteq G_\C$ with $G_\C$ simply connected, enters is the fact 
that $Z_H(\a)=Z_K(\a)$ and $H=Z_H(\a )H_0$. 
 
\ssk Recall from Subsection 2.7 the mapping 
$j(\lambda)\: \C^{\left| {\cal W}/ {\cal W}_0\right|}\to ({\cal H}_\lambda^{-\infty})^H$ 
and the intertwining  matrix $A_H(\lambda, w\lambda)^*\: ({\cal H}_{w\lambda}^{-\infty})^H\to 
({\cal H}_\lambda^{-\infty})^H$ both defined for generic $\lambda\in \a_\C^*$, and all $w\in {\cal W}$. 
For generic $\lambda$ we  define 
$$j^0(\lambda)\:  \C^{\left| {\cal W}/ {\cal W}_0\right|}\to ({\cal H}_\lambda^{-\infty})^H$$ 
by 
$$j^0(\lambda)\:=[A_H(w_0\lambda, \lambda)^*]^{-1}\circ  j(w_0\lambda)\leqno(3.3.1)$$ 
with $w_0\in {\cal W}$ the longest element. Then $j^0$ has no poles on $i\a^*$ 
(cf.\ [vdBS97b, Th.\ 1]). 
Denote by $(\e_w)_{w\in {\cal W}/ {\cal W}_0}$ the canonical basis of the Hilbert space 
$\C^{\left| {\cal W}/ {\cal W}_0\right|}$. For $w\in {\cal W}/{\cal W}_0$ define 
$$\eta_{\lambda,w}^0\:=j^0(\lambda) \e_w\in ({\cal H}_\lambda^{-\infty})^H\ .$$ 
Define a Hilbert space structure on $\Hom (\C^{\left| {\cal W}/ {\cal W}_0\right|}, {\cal H}_\lambda)$ 
using the identification 
$$\Hom (\C^{\left| {\cal W}/ {\cal W}_0\right|}, {\cal H}_\lambda)\simeq {\cal H}_\lambda\otimes 
[\C^{\left| {\cal W}/ {\cal W}_0\right|}]^*\ .$$ 
\par Write ${\cal S}(G/H)$ for the Schwartz space on $G/H$ 
and $p_{\rm mc}\: L^2(G/H)\to L^2(G/H)_{\rm mc}$ for  
the orthogonal projection on the most continuous spectrum.  
Set ${\cal S}_{\rm mc}(G/H)=p_{\rm mc}({\cal S}(G/H))$.  
Then for functions $f\in {\cal S}_{\rm mc}(G/H)$  
the Fourier transform is defined by  
$${\cal F}(f)=\left(\pi_\lambda(f)j^0(\lambda)\right )_\lambda \ .\leqno(3.3.2) $$ 
By [D98, Th.\ 3] or [vdBS97a,Cor. 18.2 and Prop. 18.3], 
${\cal F}$ extends to a $G$-equivariant unitary  
isomorphism  
$${\cal F}\: L^2(G/H)_{\rm mc} \to \int_{i\a_+^*}^\oplus 
\Hom (\C^{\left| {\cal W}/ {\cal W}_0\right|}, {\cal H}_\lambda) \ d\lambda\ . \leqno(3.3.3)$$ 
In particular, we have (using suitable normalization of measures) that 
$$\|f\|^2= \int_{i\a_+^*} 
\|{\cal F}(f)(\lambda)\|^2 \ d\lambda \leqno(3.3.4) $$ 
for all $f\in {\cal S}_{\rm mc}(G/H)$.  
\par Next we wish to describe ${\cal F}^{-1}$.  
Let $(\e_w^*)_{w\in {\cal W}/ {\cal W}_0}$ be the 
dual basis of $(\e_w)_{w\in {\cal W}/ {\cal W}_0}$. Then a section $s$ of 
$\int_{i\a_+^*}^\oplus 
\Hom (\C^{\left| {\cal W}/ {\cal W}_0\right|}, {\cal H}_\lambda) \ d\lambda$ can be written as 
$s=(\sum_{w\in {\cal W}/ {\cal W}_0} s_{\lambda, w}\otimes \e_w^*)_\lambda$ 
with $s_{\lambda, w}\in {\cal H}_\lambda$ for all $\lambda\in i\a_+^*$ and 
$w\in {\cal W}/ {\cal W}_0$. Recall that if  
$s=(\sum_{w\in {\cal W}/ {\cal W}_0} s_{\lambda, w}\otimes \e_w^*)_\lambda$ is a smooth vector, then  
$s_{\lambda, w}$ is a smooth vector in ${\cal H}_\lambda$ 
for almost all $\lambda$.  
In the sequel we will use the identification  
${\cal H}_\lambda^\infty=C^\infty(M\bs K)$. Define a 
subspace of $\left(\int_{i\a_+^*}^\oplus 
\Hom (\C^{\left| {\cal W}/ {\cal W}_0\right|}, {\cal H}_\lambda) \ d\lambda\right)^\infty$ 
by 
$$ 
{\cal H}_0=\left\{ s
\in \left(\int_{i\a_+^*}^\oplus 
\Hom (\C^{\left| {\cal W}/ {\cal W}_0\right|}, 
{\cal H}_\lambda) \ d\lambda\ \right)^\infty\quad :\quad 
\eqalign{& s \ \hbox{is $K$-finite,}\  
\supp (s) \ \hbox{is compact,}\cr 
& s\: i\a_+^*\to \Hom(\C^{\left| {\cal W}/ {\cal W}_0\right|}, 
 C^\infty(M\bs K)) \cr 
 & \hbox{is weakly smooth}. \cr}\right\}$$ 
It is not hard to see that ${\cal H}_0$ is  
a dense subspace in $\int_{i\a_+^*}^\oplus 
\Hom (\C^{\left| {\cal W}/ {\cal W}_0\right|},  
{\cal H}_\lambda) \ d\lambda$.  
Then for an element $s\in {\cal H}_0$ the inverse  
Fourier-transform is given by [D98, Th.\ 3] 
 
$${\cal F}^{-1}(s)(gH)=\int_{i\a_+^*}\sum_{w\in {\cal W}/ {\cal W}_0} \la \pi_\lambda(g^{-1}) s_{\lambda, w}, 
\eta_{\lambda, w}^0\ra \ d\lambda\ .\leqno(3.3.5)$$ 
Moreover [D98, Th.\ 3] implies that 
$${\cal F}^{-1}({\cal H}_0)\subeq  
{\cal S}_{\rm mc}(G/H)\ ;  
\leqno(3.3.6)$$ 
in particular ${\cal F}({\cal F}^{-1}(s))$ is  
given by the formula (3.3.2) for $s\in {\cal H}_0$.  
  
\Remark 3.3.1. We have normalized the invariant 
measure on $G/H$ and the measure $d\lambda$ on ${\frak a}^*$ so that (3.3.3) and 
(3.3.4) holds without any additional constants. This is possible, because we 
are only working with the principial series of representations 
and the most continuous part of the spectrum. In general, one 
has to take into account the order of several Weyl groups. 
We refer to Theorem 31 and Remark 32 in [vdB00] for general 
discussion on the normalization of measures.\qed 
 
\subheadline{3.4. Isometry of the boundary value mapping} 
 
\noindent 
In this subsection we  complete our 
discussion of the boundary value mapping begun in Subsection 3.2.

\Theorem 3.4.1. {\rm \bf(Isometry of the boundary value mapping)}
The boundary value mapping, initially defined 
by 
$$b_0\: {\cal H}^2(\Xi_H)_0\to C(G/H), \ \ b_0(f)(gH)=\lim_{t\nearrow 1} 
f(ga_tK_\C)$$ 
(cf.\ {\rm (3.2.2-3)}) extends to a $G$-equivariant isometric embedding 
$$b\:  {\cal H}^2(\Xi_H)\to L^2(G/H)_{\rm mc}\ .$$ 
 
\Proof. For each $\lambda\in i\a_+^*$ define a vector 
$\b(\lambda)\in (\C^{\left |{\cal W} 
/ {\cal W}_0\right|})^*$ by 
$$\b(\lambda)= \c(w_0\lambda)
\sum_{w\in {\cal W}/ {\cal W}_0} z_{\hbox{\fiverm H}}^{-w^{-1}(w_0\lambda+\rho)}\e_w^*\ .$$
Notice, that for $\lambda\in i\a^*$ we have
$$|z_{\hbox{\fiverm H}}^{-2w^{-1}\lambda}|=
z_{\hbox{\fiverm H}}^{-w^{-1}(\lambda+\rho)}\overline{z_{\hbox{\fiverm H}}^{-w^{-1}(\lambda+\rho)}}
=z_{\hbox{\fiverm H}}^{-w^{-1}(\lambda+\rho)}z_{\hbox{\fiverm H}}^{w^{-1}(\bar{\lambda}+\rho)}\, .$$
Therefore,
employing the Maass-Selberg relation for $\c(\lambda)$ we obtain
$$\|\b(\lambda)\|^2=|\c(\lambda)|^2 \cdot \sum_{w\in {\cal W}/ {\cal W}_0}
|z_{\hbox{\fiverm H}}^{-w^{-1}\lambda}|^2 =|\c_{G/H}(\lambda)|^2\ .$$
In particular we see that we have an  
$G$-equivariant isometric embedding of  
direct integrals  
$$\iota\: \int_{i\a_+^*}^\oplus {\cal H}_\lambda\ d\mu(\lambda) 
\to  \int_{i\a_+^*}^\oplus \Hom(\C^{\left |{\cal W} 
/ {\cal W}_0\right|}, {\cal H}_\lambda)\ d\lambda,  
\ \  s=(s_\lambda)_\lambda\mapsto (s_\lambda\otimes \b(\lambda))_\lambda\ .$$ 
{}From the definition of the spaces ${\cal H}^2(\Xi_H)_0$ 
and ${\cal H}_0$ it is then clear that
$$\iota \left(\Phi^{-1}({\cal H}^2(\Xi_H)_0)\right)\subeq 
{\cal H}_0\ .\leqno(3.4.1)$$ 

\par As
$$\iota (s)_\lambda = \c (w_0\lambda)\sum_{w\in {\cal W}/{\cal W}_0}
z_{\hbox{\fiverm H}}^{-w^{-1}(w_0\lambda+\rho)}s_\lambda\otimes \e_w^*\ , $$
we get by (3.2.3), (3.3.5),  Theorem
2.6.1 and Theorem 2.7.1 that
$$\eqalign{[{\cal F}^{-1}(\iota (s))](gH) &=\int_{i\a_+^*}
\sum_{w\in {\cal W}/{\cal W}_0}
\c (w_0\lambda)
z_{\hbox{\fiverm H}}^{-w^{-1}(w_0\lambda+\rho)}
\la \pi_\lambda (g^{-1})s_\lambda, \eta^0_{\lambda, w}\ra \ d\lambda\cr
&=\int_{i\a_+^*}
\sum_{w\in {\cal W}/{\cal W}_0}
\la \pi_\lambda (g^{-1})s_\lambda, \c (w_0\bar{\lambda})
z_{\hbox{\fiverm H}}^{w^{-1}(w_0\bar{\lambda}+\rho)}\eta^0_{\lambda, w}\ra \ d\lambda\cr
&=\int_{i\a_+^*}
\la \pi_\lambda (g^{-1})s_\lambda, v_{{\hbox{\fiverm H}, \scriptscriptstyle{\lambda}}}\ra \ d\lambda
\cr
&=b_0(s)(gH) }
$$

\par From this and (3.3.6) it follows that
$b_0(f)\in S_{\rm mc}(G/H)$; in particular we
have
$b_0(f)\in L^2(G/H)_{\rm mc}$. Finally,
$$\|b_0(f)\|_{L^2(G/H)_{\rm mc}}=
\|{\cal F}^{-1} (\iota (s))\|=\|\iota(s)\|=\|s\|=
\|\Phi(s)\|=\|f\|$$ 
as ${\cal F}$, $\iota$ and $\Phi$ are isometric.  
This completes the proof of the theorem.\qed

\Remark 3.4.2. The domain $\Xi_H$ is maximal in the sense that 
generic functions in ${\cal H}^2(\Xi_H)$ do not extend holomorphically 
over $\Xi_H$.\qed

\subheadline{3.5. Concluding remarks and the example of $G=\Sl(2,\R)$} 
 
\noindent 
In [GK\'O01] we defined a Hardy space ${\cal H}^2(\Xi)$ on $\Xi$ for the  
cases where $\Xi=\Xi_H$.  Let us briefly summarize its construction in order  
to put it into perspective with the results in this section. 
 
\par Geometrically the situation $\Xi=\Xi_H$ is equivalent to the fact  
that $\Xi$ is homogeneous 
for a bigger Hermitian group $S\supeq G$ (cf.\ [KS01b]). More precisely, if $U<S$ denotes an appropriate  
maximal compact subgroup with $K\subeq U$ then $\Xi$ is $G$-biholomorphic  
to the Hermitian symmetric space $S/U$. For example if $G$ is Hermitian, then  
$S=G\times G$ and $\Xi\simeq G/K\times \oline{G/K}$.  
 
\par The assumption $\Xi=\Xi_H$ thus allows us to identify $\Xi$ with a  bounded  
symmetric domain ${\cal D}\simeq S/U$. Within this identification one shows that  
$\partial_d\Xi\simeq G/H$ becomes a Zariski-open subset in the Shilov boundary  
$\partial_s{\cal D}$ of ${\cal D}$.

\par The identification of $\Xi$ with ${\cal D}$ was used in [GK\'O01] in a crucial way:  
One can transfer the action of an appropriate compression-semigroup $\Gamma\supeq G$ on ${\cal D}$ 
to $\Xi$ and use this to give a definition of a Hardy space as follows: 
 
$${\cal H}^2(\Xi)=\{ f\in {\cal O}(\Xi)\: \|f\|^2=\sup_{\gamma\in \Int\Gamma} 
\int_{G/H} |f(\gamma gz_{\hbox{\fiverm H}})|^2\ dgH<\infty\}\ .\leqno(3.5.1)$$ 
In [GK\'O01] we have shown -- with entirely different methods -- that the Hardy space defined as  
in (3.5.1) has the following properties: 
\msk 
\item{(3.5.2)} ${\cal H}^2(\Xi)$ is a Hilbert space of holomorphic functions 
featuring an isometric boundary value mapping $b\:{\cal H}^2(\Xi)\into L^2(G/H)_{\rm mc}$.  
Moreover, $\im b$ is a multiplicity one subspace of {\it full spectrum}.  
\item{(3.5.3)} The Hardy space ${\cal H}^2(\Xi)$ is $G$-isometric to  
$L^2(G/K)$ through a transform of Segal-Barg\-mann type.   
\item{(3.5.4)} ${\cal H}^2(\Xi)$ is $G$-isometric to the classical Hardy space  
${\cal H}^2({\cal D})$ through an explicitely given mapping.  
\msk 
In particular for $\Xi=\Xi_H$ it follows from Theorem 3.4.1 and  
(3.5.2) that the definition of (3.5.1) coincides 
with our spectral definition of the Hardy space in Definition 3.2.2.  For the  
cases where $\Xi\neq \Xi_H$ there is no apparent semigroup action on $\Xi_H$ and a definition  
of ${\cal H}^2(\Xi_H)$ in the flavour of (3.5.1) seems presently not possible.  
 
\par Notice that (3.5.3) implies that the Plancherel measure of ${\cal H}^2(\Xi)$ has  
support equal to $i\a_+^*$. However, in [GK\'O01] we could not determine this  
measure explicitely. With the new approach given in this section this difficulty  
is already taken care of with the definition of the Hardy space.  
\par The explicit isomorphism of ${\cal H}^2(\Xi)$ with the classical  
Hardy space ${\cal H}^2({\cal D})$ allows us to find also a nice closed  
expression for the Cauchy-Szeg\"o function $\Psi$ (cf.\ [GK\'O01, Th.\ 5.7 and  
Ex.\  5.10]). Combining this closed expression with the spectral resolution  
of $\Psi$ in Theorem 3.2.4 one obtains interesting identities for  
(generalized) hypergeometric functions. For example for $G=\Sl(2,\R)$ one obtains 
the following formula: 
$${1-\tanh^2 t\over 1+\tanh^2 t} ={\pi\over 2}\int_0^\infty F\left({1\over 4}+i{\lambda\over 4},  
{1\over 4}-i{\lambda\over 4}, 1; -\sinh^2\left(2t+i{\pi\over 2}\right)\right) 
\cdot\left|{\Gamma\left({i\lambda+1\over2}\right)\over \Gamma\left({i\lambda\over2}\right)} \right|^2 
\ {d\lambda\over \cosh {\pi\over 2}\lambda}$$ 
for all $t\in \R+i]-{\pi\over 4}, {\pi\over 4}[$. Here $F$ denotes the Gau\3 hypergeometric function.

\sectionheadline{A. Appendix: Analytic vectors for representations} 
 
In this appendix  we will summarize  some facts on analytic vectors  
for representations. None of the results collected below is new,  
however some of them might be hard to find explicitely in the  
literature. In order to keep the exposition short, we will omit  
proofs and often do not make the most general assumptions.  
A more detailed account containing complete proofs can be found  
in the forthcoming survey [K\'O03].   
 
\subheadline{A.1. Definition and topology of analytic vectors} 
 
Throughout this appendix $G$ will denote a connected unimodular Lie group  
with $G\subeq G_\C$.  
 
\par Let $E$ be a complex Banach space and  
$\Gl(E)$ the group of continuous invertible operators on $E$.  
By a (Banach) representation of $(\pi, E)$ of $G$ we will understand  
a group homomorphism $\pi\: G\to \Gl(E)$ such that  
for all $v\in E$ the orbit mapping  
$$\gamma_v\: G\to E, \ \ g\mapsto \pi(g)v$$ 
is continuous.  
 
\par A vector $v\in E$ is called {\it analytic} if $\gamma_v$ is an  
analytic $E$-valued map or, equivalently, if there exists  
an open neighborhood $U$ of $\1$ in $G_\C$ and a $G$-equivariant  
holomorphic mapping  
$$\gamma_{v,U}\: GU\to E$$ 
such that $\gamma_{v,U}(\1)=v$. In particular, $\gamma_{v, U}\res_G=\gamma_v$.  
 
\par The vector space of all analytic vectors  for $(\pi, E)$ is denoted  
by $E^\omega$. We recall a fundamental result of Nelson which  
states that $E^\omega$ is dense in $E$.  
 
\ssk  Next we are going to recall the definition of the topology  
on $E^\omega$.  
\par For a complex manifold 
$M$ let us denote by ${\cal O}(M,E)$ the space of all $E$-valued  
holomorphic mappings on $E$. Topologically we consider ${\cal O}(M,E)$  
as a Fr\'echet space with the topology of compact convergence. 
 
\par  For any open neighborhood $U$  
of $\1$ in $G_\C$ we write $E_U$ for the subspace  
of $E^\omega$ for which $\gamma_{v,U}$ exists. Then we obtain a linear  
embedding  
 
$$\eta_U\: E_U\to {\cal O}(GU, E), \ \ v\mapsto \gamma_{v,U}\ .$$ 
The image of $\eta_U$ is closed and hence ${\cal O}(GU,E)$ induces  
a Fr\'echet topology on $E_U$. Notice that for $U_1\subeq U_2$  
we obtain a continuous embedding $E_{U_2}\to E_{U_1}$ via restriction.  
Thus  
$$E^\omega=\lim_{U\to \{1\}} E_U=\bigcup_U E_U $$  
and we can equip $E^\omega$ with the inductive limit topology, i.e. 
the finest topology on $E^\omega$ for which all inclusion  
mappings $E_U\to E^\omega$ become continuous. Notice that this  
turns $E^\omega$ into a locally convex topological vector space.   
 
\par By $E^{-\omega}$ we will denote the antidual of $E^\omega$, i.e. 
the space of all antilinear continuous functionals on $E^\omega$.  
The space $E^{-\omega}$ is referred to as the space of {\it hyperfunction  
vectors} of the representation $(\pi, E)$. We equip  
$E^{-\omega}$ with the topology of bounded convergence.

\subheadline{A.2. Analytic vectors for $L^1(G/H)$} 
 
\par Let $H<G$ be a closed subgroup such that $G/H$  
carries a $G$-invariant measure. We write $L^1(G/H)$ for the corresponding  
Banach space of integrable functions and $(L, L^1(G/H))$ for the left regular representation  
of $G$ on $L^1(G/H)$, i.e.,  
$$(L(g)f)(xH)=f(g^{-1}xH)\qquad (g,x\in G, f\in L^1(G/H))\ .$$ 
 
Further it is convenient to assume that $G/H\subeq G_\C/ H_\C$.  
Then we have the following characterization of the analytic 
vectors:  
 
\Proposition A.2.1. Let $U$ be an open neighborhood of $\1$ in $G_\C$.  
Then $f\in L^1(G/H)_U$ if and only if  
there exists a holomorphic function $\tilde f$ on the open set  
$$U^{-1} GH_\C/ H_\C \subeq G_\C/ H_\C$$  
with the following properties: 
\item{(1)} $\tilde f\res_{G/H}=f$.  
\item{(2)} For all $x\in U$ the map  
$$\tilde f_x\: G/H\to\C, \ \ gH\mapsto \tilde f(x^{-1}gH)$$ 
belongs to $L^1(G/H)$.  
\item{(3)} For all compact subsets $U^c\subeq U$ we have  
$$\sup_{x\in U^c} \|\tilde f_x\|<\infty\ .$$ \qed  
 
There are two types  of homogeneous spaces  $G/H$ which will 
be of particular interest for us. The first is when $H={\1}$. Then $L^1(G)^\omega$ 
denotes the analytic for the left regular representation of $G$ on $L^1(G)$.  
The second case is for $G=H\times H$ and $H<G$ the diagonal subgroup. In this  
case $G/H\simeq H$ and $L$ becomes left-right regular representation of $H\times H$  
on $H$. Here we shall write $L^1(H)^{\omega, \omega}$ for the analytic vectors.

\subheadline{A.3. Averaging properties} 
 
Recall that the average map  
$$C_c(G)\to C_c(G/H), \ \ f\mapsto f^H; f^H(xH)=\int_H f(xh)\ dh$$ 
is contiunuous and onto. Further, this map extends to a surjective contraction of Banach spaces  
$L^1(G)\to L^1(G/H)$. We will show that the  averaging operator  
maps analytic vectors into analytic vectors. 
 
\par A standard application of the Bergman estimate gives: 
 
\Lemma A.3.1. Let $U\subeq G_\C$ be an open neighborhood  
of $\1$. Then for any pair of compact subsets 
$U_1, U_2\subeq U$ with $U_1\subeq \Int U_2$ there  
exists a constant $C>0$ such that for all $f\in L^1(G)_U$ we have that  
$$(\forall x\in U_1^{-1}G) \qquad \int_H |\tilde f(xh)|\ dh\leq  
C \sup_{x\in U_2}\|\tilde f_x\|,  $$ 
where $\tilde f$ denotes the extension of $f$ to a holomorphic  
function on $U^{-1}G$ (cf.\ Proposition A.2.1).\qed  
 
Combining Lemma A.3.1 with Proposition A.2.1 we obtain:  
 
\Proposition A.3.2.  Let $U\subeq G_\C$ be an open neighborhood  
of $\1$.  Then for every  
$f\in L^1(G)_U$ and $g\in G$ the integral 
$f^H(g)=\int_H f(gh)\ dh$ converges absolutely and $f^H\in L^1(G/H)_U$.  
In particular, there is a well defined mapping  
$$L^1(G)^\omega\to L^1(G/H)^\omega, \ \ f\mapsto f^H\ .$$\qed

\subheadline{A.4. Mollifying properties} 
 
In this section $E={\cal H}$ will be a Hilbert space and $(\pi, {\cal H})$ a unitary 
representation of $G$.  
For $f\in L^1(G)$ one defines a continuous operator $\pi(f)\: {\cal H}\to {\cal H}$ by  
 
$$\pi(f)v=\int_G f(g) \pi(g)v\ dg\qquad (v\in {\cal H})\ .$$ 
Notice that this defines a $*$-representation of the Banach algebra 
$L^1(G)$, i.e. we have $\pi(f*g)=\pi(f)\pi(g)$ and $\pi(f)^*=\pi(f^*)$ with  
$f^*(x)=\oline {f(x^{-1}) }$.  
 
\par Recall that $L^1(G)^{\omega, \omega}$ denotes the analytic  
vectors for the left-right regular representation of $G\times G$ on  
$L^1(G)$. It is easy to see  
that $L^1(G)^{\omega,\omega}$ is $*$-closed subalgebra of $L^1(G)$.  
 
\par Let $f\in L^1(G)^{\omega,\omega}$. It follows readily  
from Proposition A.2.1 and the definition of analytic vectors 
that $\pi^\omega(f)$ maps ${\cal H}$ continuously into ${\cal H}^\omega$.  
In particular the restriction $\pi^\omega(f)\:=\pi(f)\res_{{\cal H}^\omega}$ 
gives rise to a continuos operator $\pi^\omega(f)\: {\cal H}^\omega\to {\cal H}^\omega$. 
Hence we have an algebra representation: 
 
$$\pi^\omega\: L^1(G)^{\omega,\omega}\to \End({\cal H}^\omega), \ \ f\mapsto \pi^\omega(f)\ .$$  
The corresponding dual representation is given by  
 
$$\pi^{-\omega}\: L^1(G)^{\omega,\omega}\to \End({\cal H}^{-\omega}); \ \   
\pi^{-\omega}(f)\lambda=\lambda\circ \pi^\omega(f^*)\ .$$ 
Another application of Proposition A.2.1 then gives us  
the mollifying property:  
 
\Proposition  A.4.1. Let $(\pi, {\cal H})$ be a unitary  
representation of a unimodular Lie group $G$. Then we have for all  
$f\in L^1(G)^{\omega,\omega}$ that   
$$\pi^{-\omega}(f){\cal H}^{-\omega}\subeq {\cal H}^\omega\ .$$\qed  
 
\nin {\bf Note:} For $f\in L^1(G)^{\omega, \omega}$, it is often convenient to write $\pi(f)$  
instead of $\pi^{-\omega}(f)$. We will use this convention throughout Section 2 in the main text.

\def\entries{

\[AG90 Akhiezer, D.\ N., and S.\ G.\ Gindikin, {\it On Stein 
extensions of 
real symmetric spaces}, 
Math.\ Ann.\ {\bf 286}, 1--12, 1990 
 
\[vdB88 van den Ban, E., {\it The principal series for a reductive symmetric space I, 
$H$-fixed distribution vextors}, Ann. sci. \'Ec. Norm. Sup. {\bf 4}, {\bf 21} 
(1988), 359--412 
 
\[vdB00 ---, {\it The Plancherel theorem for a reductive symmetric space}, Lectures 
for the European School of Group Theory. August 14--26,  
2000, SDU-Odense University. 
http://www.math.uu.nl/people/ban/publ.html 
 
\[vdBD88 van den Ban, E., and P. Delorme, {\it Quelques propri\'et\'es des repr\'esentations  
sph\'eriques pour les espaces sym\'etriques r\'eductifs}, J. Funct. Anal. {\bf 80} (1988),  
284--307 
 
\[vdBS97a van den Ban, E., and H.\ Schlichtkrull, {\it The most 
continuous part of the Plan\-che\-rel decomposition for a reductive 
symmetric space}, Ann. of Math. {\bf (2) 145} (1997), no. {\bf 2}, 
267--364 
 
\[vdBS97b ---, {\it Fourier transform on a semisimple symmetric space},  Invent. Math. {\bf 130} (1997), no. {\bf 3}, 
517--574 
 
\[BD92 Brylinski, J.-L., and P. Delorme, {\it Vecteurs distributions $H$-invariants pour les s\'eries 
principales g\'en\'eralis\'ees d'espaces sym\'etriques r\'eductifs et 
prolongement m\'eromorphe d'int\'egralesd'Eisenstein}, Invent. Math. {\bf 109} (1992), no. {\bf 3}, 619--664 
 
\[D98  Delorme, P., {\it Formule de Plancherel pour les espaces 
sym\'etriques r\'eductifs},  Ann. of Math. {\bf (2) 147} 
(1998), no. {\bf 2}, 417--452 
 
\[FT99 Faraut, J., and E.  G.  F.  Thomas, {\it Invariant Hilbert spaces of 
holomorphic functions}, J.  Lie Theory {\bf 9} (1999), no.  {\bf 2}, 383--402

\[GK02a Gindikin, S., and B.\ Kr\"otz, {\it Complex crowns of Riemannian 
symmetric spaces and non-compactly causal symmetric spaces}, 
Trans. Amer. Math. Soc. {\bf 354} (2002), no. {\bf 8}, 
3299--3327 
 
\[GK02b ---, {\it Invariant Stein domains in Stein symmetric spaces 
and a non-linear complex convexity theorem}, IMRN {\bf 18} (2002), 959--971 
 
\[GK\'O083 Gindikin, S., B.\ Kr\"otz and G.\ \'Olafsson, {\it Hardy spaces for 
non-compactly causal symmetric spaces and the most continuous spectrum}, Math. Ann. {\bf 327} (2003), 25--66
 
\[H84 Helgason, S., ``Groups and Geometric Analysis'', Academic Press, 1984  
 
\[H\'O96 Hilgert, J.\ and 
G.\ \'Olafsson, ``Causal Symmetric Spaces, Geometry and 
Harmonic Analysis,'' Acad. Press, 1996 
 
\[K99 Kr\"otz, B., {\it The Plancherel theorem for biinvariant Hilbert 
spaces}, Publ.  Res.  Inst.  Math.  Sci.  {\bf 35} (1999), no.  {\bf 1}, 
91--122 
 
\[K\'O03 Kr\"otz, B., and G. \'Olafsson, {\it Analytic vectors for representations -- a  
survey}, in preparation 
 
\[KS01a Kr\"otz, B., and R.J. Stanton, {\it Holomorphic extension of 
representations: (I) 
automorphic functions}, Annals of Mathematics, to appear 
 
\[KS01b Kr\"otz, B., and R.J. Stanton, {\it Holomorphic extensions of 
representations: (II) geometry and harmonic analysis}, preprint 
 
\[M79 Matsuki, T., {\it The orbits of affine symmetric spaces under the action of minimal parabolic 
subgroups}, J. Math. Soc. Japan {\bf 31} (1979), 331--357 
 
\[\'O87 \'Olafsson, G., {\it Fourier and Poisson transformation associated to 
semisimple symmetric space}, Invent Math. {\bf 90} (1987), 605--629

}

{\sectionheadline{\bf References} 
\frenchspacing 
\entries\par} 
\dlastpage  
\bye